\documentclass[11pt]{amsart}

\usepackage{a4wide, amsmath, amsfonts, amssymb, mathrsfs, amsthm, bbm, stmaryrd}
\usepackage[hyperfootnotes=false]{hyperref}

\allowdisplaybreaks[2]

\numberwithin{equation}{section}

\newtheorem{theorem}{Theorem}[section]
\newtheorem{lemma}[theorem]{Lemma}

\newtheorem{proposition}[theorem]{Proposition}
\newtheorem{definition}[theorem]{Definition}

\theoremstyle{remark}
\newtheorem{remark}[theorem]{Remark}

\renewcommand{\P}{\mathbb{P}}
\renewcommand{\d}{\mathrm{d}}
\renewcommand{\epsilon}{\varepsilon}
\newcommand{\dd}{\,\mathrm{d}}
\newcommand{\R}{\mathbb{R}}

\newcommand{\N}{\mathbb{N}}

\newcommand{\E}{\mathbb{E}}
\newcommand{\1}{\mathbf{1}}
\newcommand{\D}{\mathrm{D}}
\newcommand{\sign}{\textup{sign}}

\newcommand\bL{{\mathcal K}}
\newcommand\cS{\mathcal{S}}
\global\long\def\rbr#1{\left(#1\right)}
\global\long\def\sbr#1{\left[#1\right]}
\global\long\def\cbr#1{\left\{#1\right\}}
\global\long\def\TTV#1#2#3{\text{TV}^{#3}\!\rbr{#1,#2}}

\global\long\def\Ucross#1#2#3{\text{u}^{#1}\!\rbr{#2,#3}}
\global\long\def\Dcross#1#2#3{\text{d}^{#1}\!\rbr{#2,#3}}
\global\long\def\cross#1#2#3{\text{n}^{#1}\!\rbr{#2,#3}}
\global\long\def\Cross#1#2#3{{\emph{\text{n}}^{#1}\!\rbr{#2,#3}}}
\global\long\def\Ucrosstilde#1#2#3{\text{\~{u}}^{#1}\!\rbr{#2,#3}}
\global\long\def\Dcrosstilde#1#2#3{\text{\~{d}}^{#1}\!\rbr{#2,#3}}
\global\long\def\crosstilde#1#2#3{\text{\~{n}}^{#1}\!\rbr{#2,#3}}

\title[Local times and Tanaka--Meyer formulae for c{\`a}dl{\`a}g paths]{Local times and Tanaka--Meyer formulae \\for c{\`a}dl{\`a}g paths}

\author[\L ochowski]{Rafa\l{} M. {}\L ochowski}
\address{Rafa\l{} M. {\L}ochowski, Warsaw School of Economics, Poland}
\email{rlocho314@gmail.com}

\author[Ob{\l}{\'o}j]{Jan Ob{\l}{\'o}j}
\address{Jan Ob{\l}{\'o}j, University of Oxford, United Kingdom}
\email{obloj@maths.ox.ac.uk}

\author[Pr{\"o}mel]{David J. Pr{\"o}mel}
\address{David J. Pr{\"o}mel, University of Mannheim, Germany}
\email{proemel@uni-mannheim.de}

\author[Siorpaes]{Pietro Siorpaes}
\address{Pietro Siorpaes, Imperial College London, United Kingdom}
\email{p.siorpaes@imperial.ac.uk}

\date{\today.}

\begin{document}

\begin{abstract}
  Three concepts of local times for deterministic c{\`a}dl{\`a}g paths are developed and the corresponding pathwise Tanaka--Meyer formulae are provided. For semimartingales, it is shown that their sample paths a.s. satisfy all three pathwise definitions of local times and that all coincide with the classical semimartingale local time. In particular, this demonstrates that each definition constitutes a legit pathwise counterpart of probabilistic local times. The last pathwise construction presented in the paper expresses local times in terms of normalized numbers of interval crossings and does not depend on the choice of the sequence of grids. This is a new result also for c{\`a}dl{\`a}g semimartingales, which may be related to previous results of Nicole El~Karoui~\cite{ElKaroui1978} and Marc Lemieux~\cite{Lemieux1983}. 
\end{abstract}

\maketitle



\section{Introduction}

Stochastic calculus, with its foundational notions developed by Kyiosi It\^o in the 1940s, is a \emph{par excellence} probabilistic endeavour. The stochastic integral, the integration by parts formula -- these basic building blocks are to be understood almost surely, and so is the edifice they span. This thinking has proved to be exceedingly powerful and fruitful, and underpins many beautiful developments in probability theory since then. Nevertheless, for decades now, mathematicians have been trying to develop a more analytic, pathwise understanding of these probabilistic objects. On one hand, this was, and is, driven by mathematical curiosity. The classical calculus remains an irresistible reference point and, e.g., in developing a notion of an integral it is important to understand when and how it can be seen as a limit of its Riemann sums. On the other hand, this was, and is, driven by applications. Stochastic differential equations have became a ubiquitous tool for mathematical modelling from physics, through biology to finance. Yet, they do not offer the same level of path-by-path description of the system's evolution as the classical differential equations do. This becomes particularly problematic if one needs to work simultaneously with many probability measures, possibly mutually singular. One field where this proves important, and which has driven renewed interest in pathwise stochastic calculus, is robust mathematical finance, see for example \cite{Davis2014} and the references therein. Both of the above reasons -- mathematical curiosity and possible applications -- are important for us. We add to this literature and develop a pathwise approach to stochastic calculus for c{\`a}dl{\`a}g paths using local times.

In his seminal paper \cite{Follmer1981}, F{\"o}llmer introduced, for twice continuously differentiable $f\colon \R\to \R$, a non-probabilistic version of the It{\^o} formula  
\begin{align*}
  f(x_t) -f(x_0)  = \int_0^t  f^\prime(x_{s-}) \dd x_s + \frac{1}{2} \int_0^t f^{\prime \prime} (x_s) \dd [x]^c_s 
  + J_t^f(x), \quad t\in [0,T], 
\end{align*} 
where $x\colon [0,T]\to \R$ is c{\`a}dl{\`a}g and possesses a suitably defined quadratic variation $[x]$ such that, for $0\leq t\leq T$,
\begin{equation*}
  [x]_t=[x]^c_t+\sum_{0<s\le t}(\Delta x_s)^2 , \text{ where } \Delta x_t := x_{t}-x_{t-} ,
\end{equation*}
and $J_t^f(x)$ is defined by the following absolutely convergent series
\begin{equation*}
  J_t^f(x):= \sum_{0<s\le t}\big ( \Delta f(x_{s})-f^\prime (x_{s-} )\Delta x_s \big) .
\end{equation*}
In particular, this leads to a pathwise definition of the ``stochastic'' integral $\int_0^t f^\prime(x_{s-}) \dd x_s$, assuming $[x]$ exists. Soon after, Stricker \cite{Str81}, showed that one could not extend the above to all continuous functions $f$. This could only be done adopting a much more bespoke discretisation and probabilistic methods, see for example \cite{Bichteler1981,Karandikar1995}. Accordingly, the main remaining challenge was to understand the case of functions $f$ which are not twice continuously differentiable but are weakly differentiable, in some sense. In probabilistic terms, this realm is covered by the Tanaka--Meyer formula. 

For continuous paths F{\"o}llmer's pathwise It{\^o} formula was generalized to a pathwise Tanaka--Meyer formulae in the early work of \cite{Wuermli1980} and more recently in \cite{Perkowski2015} and in \cite{Davis2018}, who offered a comprehensive study. Furthermore, we refer to \cite{Geman1980} and \cite{Bertoin2014,Davis2014} for related work in a pathwise spirit. Our contribution here is to study this problem for c{\`a}dl{\`a}g paths. Jump processes, e.g., L{\'e}vy processes, are of both theoretical and practical importance and, as stressed above, our study is motivated by both mathematical curiosity as well as applications. Already in the classical, probabilistic, setting stochastic calculus for jump processes requires novel insights over and above the continuous case. This was also observed in recent works focusing on F{\"o}llmer's It{\^o} calculus for c{\`a}dl{\`a}g paths, see \cite{Chiu2018} and \cite{Hirai2019}. We face the same difficulty, which of course makes our study all the more interesting. In particular, we need more information and new ideas to handle jumps. This is consistent with the definition of quadratic variation for c{\`a}dl{\`a}g paths, cf. \cite{Chiu2018}. 

Our non-probabilistic versions of Tanaka--Meyer formula, extend the above It{\^o} formula allowing for functions~$f$ with weaker regularity assumptions than $C^2$. More precisely, we derive pathwise formulae 
\begin{equation*}
  f(x_t)-f(x_0)  
  =\int_{0}^t f^\prime(x_{s-})\dd x_s+ \frac{1}{2} \int_{\R} L_t(x,u)  f^{''}(\mathrm{d}u)  
  +J_t^f(x),\quad t\in [0,T], 
\end{equation*}
for twice weakly differentiable functions $f$, supposing that the c{\`a}dl{\`a}g path $x$ possesses a suitable pathwise local time~$L(x)$. As in the case of the It{\^o} formula, there exists no unique pathwise sense to understand such a formula, see also Remark~\ref{rem:pathwise Ito formula} below. We develop three natural pathwise approaches to local times and, consequently, to their stochastic calculus. First, we start with the key property relating local times and quadratic variation: the time-space occupation formula, and use it to define pathwise local times. Second, in the spirit of \cite{Follmer1981,Wuermli1980}, we discretise the path along a sequence of partitions and obtain local times as limits of discrete level crossings and stochastic integrals as limits of their Riemann sums. Finally, we discretise the integrand via the Skorokhod map which provides a natural approximation of the ``stochastic'' integral and links to the concept of truncated variation. In all of the three cases we show that a pathwise variant of the Tanaka--Meyer formula holds. Further, we prove that for a c{\`a}dl{\`a}g semimartingale, all three constructions coincide a.s. with classical local times. This shows that all three approaches are legitimate extensions of the classical stochastic results to pathwise analysis. Each has its merits and limitations which we explore in detail. Our aim is to provide a comprehensive understanding of how to deal with jumps in the context of pathwise Tanaka--Meyer formulae. We thus do not seek further extensions of the setup, e.g., to cover time-dependent functions $f$, cf. \cite{Feng2006}, path-dependent functions, cf. \cite{Cont2010,Imkeller2015,Saporito2018}, nor to develop higher order local times in the spirit of \cite{Cont2019} for c{\`a}dl{\`a}g paths. These, while interesting, would distract from the main focus of the paper and are left as avenues for future research.

\medskip

\noindent{\bf Outline:} In Section~\ref{sec:tanaka and local times} we propose three notions of local times for c{\`a}dl{\`a}g paths and establish the corresponding Tanaka--Meyer formulae. Then, in Section~\ref{sec:construction of local limes}, we show that sample paths of semimartingales almost surely possess such local times and all three definitions agree a.s. in the classical stochastic world.

\medskip

\noindent{\bf Acknowledgement:} This project was generously supported by the European Research Council under (FP7/2007-2013)/ERC Grant agreement no. 335421. The research of RM{\L} was partially supported by the National Science Centre (Poland) under the grant agreements no. 2016/21/B/ST1/0148 and no. 2019/35/B/ST1/0429. JO is grateful to St John's College Oxford for their support, and to the Sydney Mathematical Research Institute, where the final stages of this research were completed, for their hospitality.

\section{Pathwise local times and Tanaka--Meyer formulae}\label{sec:tanaka and local times}

The first non-probabilistic version of It{\^o}'s formula and the corresponding notion of pathwise quadratic variation of c{\`a}dl{\`a}g paths was introduced by H. F{\"o}llmer in the seminal paper~\cite{Follmer1981}. Before providing non-probabilistic versions of Tanaka--Meyer formulae and introducing the corresponding pathwise local times, we recall in the next subsection some results from~\cite{Follmer1981}.

\subsection{Quadratic variation and the F{\"o}llmer--It{\^o} formula}

For $T\in (0,\infty)$, let $D([0,T];\R)$ be the space of all c{\`a}dl{\`a}g (RCLL) functions $x\colon[0,T]\to\R$, that is, $x$ is right-continuous and possesses finite left-limits at each $t\in [0,T]$. For $x\in D([0,T];\R)$ we set $x_{t-}:=\lim_{s<t, s\to t}x_s$ for $t\in (0,T]$, $x_{0-}:=x_0$ and $\Delta x_s := x_s-x_{s-}$ for $s\in [0,T]$.

In order to define the summation over the jumps of a c{\`a}dl{\`a}g function, we need the concept of summation over general sets, see for example \cite[p.77-78]{Kelly1975}. Let $I$ be a set, let $b\colon I\to \R$ be a real valued function and let $\mathcal{I}$ be the family of all finite subsets of $I$. Since $\mathcal{I}$ is directed when endowed with the order of inclusion~$\subseteq$, the summation over $I$ can be defined by 
\begin{equation}\label{Sum AC}
  \sum_{i\in I} b_i:=\lim_{\Gamma \in \mathcal{I}}\sum_{i\in \Gamma} b_i
\end{equation} 
as limit of a net, i.e., $\lim_{\Gamma \in \mathcal{I}}\sum_{i\in \Gamma} b_i=:l \in [-\infty,\infty]$ exists if, for any neighbourhood\footnote{The space $[-\infty,\infty]$ is given the usual topology which makes it isomorphic to $[-1,1]$; in particular one can take $(x-\epsilon, x+\epsilon)$ (resp. $(M,+\infty)$, resp. $(-\infty,-M)$), where $0<\epsilon<1<M<\infty$, as a neighbourhood basis of $x\in \R$ (resp. $+\infty$, resp. $-\infty$), and metrize this topology with the distance $d(x,y):=\arctan(|x-y|)$, where $\arctan(\pm \infty):=\pm 1$, $x,y \in [-\infty,\infty]$.} $V_l$ of $l$, there is  $\Gamma\in \mathcal{I}$ such that for all $\tilde \Gamma\in \mathcal{I}$ such that $\tilde \Gamma \geq \Gamma$ (i.e., $\tilde \Gamma \supseteq \Gamma$) one has $ \sum_{i\in \tilde\Gamma} b_i\in V_l$. If $b_i\geq 0$ for all $i\in I$, then it is easy to see that 
\begin{align}\label{Sum positive}
  \exists \sum_{i\in I} b_i=\sup \bigg\{\sum_{i\in J} b_i : J  \in \mathcal{I} \bigg \} \in [0,\infty] .
\end{align}  
We say that the series $\sum_{i\in I} b_i$ is \emph{absolutely summable} if the limit $\sum_{i\in I} |b_i|$ (which always exists, by \eqref{Sum positive}) is finite, in which case also the limit \eqref{Sum AC} exists and satisfies $|\sum_{i\in I} b_i| \leq \sum_{i\in I} |b_i| $, and there exists\footnote{Since $I_n:=\{i\in I: |b_i| \geq 1/n\}$ is finite for each $n$, because $\frac{\# I_n}{n} \leq \sum_{i\in I_n} |b_i|\leq \sum_{i\in I} |b_i|<\infty$.} a countable subset $K \subseteq I$ s.t. $b_i=0$ if $i\in I \setminus K$.

For a continuous function $f\colon \R\to\R$ possessing a left-derivative $f^\prime$, we now set
\begin{equation}\label{eq:sum of jumps}
  J_t^f(x):= \sum_{0<s\le t}\big ( \Delta f(x_{s})-f^\prime (x_{s-} )\Delta x_s \big),
\end{equation}
provided the sum exists. Furthermore, the space of continuous functions $f\colon \R\to \R$ is denoted by $C(\R):=C(\R;\R)$, the space of twice continuously differentiable functions by $C^2(\R):=C^2(\mathbb{R};\mathbb{R})$ and the space of smooth functions by $C^\infty(\R):=C^\infty(\R;\R)$.

A \textit{partition} $\pi = (t_j)_{j=0}^N $ is a finite sequence such that $0 = t_0 < t_1 < \dots <t_N=T$ (for some $N\in \N$). We write $|\pi|:= \max_{ j \in \mathbb{N}} \vert t_j - t_{j-1}\vert$ for its mesh size and define $\pi(t):=\pi\cap [0,t]$ the restriction of $\pi$ to $[0,t]$. A sequence of partitions $(\pi^n)_{n \in \mathbb{N}}$ is said to be \textit{refining} if for all $t_j \in \pi^n$ we also have $t_j \in \pi^{n+1}$ and a refining sequence $(\pi^n)_{n \in \mathbb{N}}$ is said to exhaust the jumps of $x$ if for all $t\in [0,T]$ with $\Delta x_t \neq 0$, $t\in \pi^n$ for $n$ large enough. The Dirac measure at $t \in [0,T]$ is denoted by $\delta_t$.

\begin{definition}
  Let $(\pi^n)_n$ be a sequence of partitions such that $\lim_{n\to\infty}|\pi^n|=0$. A function $x\in D([0,T];\R)$ has \emph{quadratic variation} $[x]$ along $(\pi^n)_n$ if  the sequence of discrete measures 
  \begin{equation*}
    \mu_n := \sum_{t_j \in \pi^n} (x_{t_{j+1}} - x_{t_j})^2 \delta_{t_j}
  \end{equation*}
  converges weakly\footnote{Meaning that $\int_0^T h \dd \mu_n\to \int_0^T h \dd \mu$ for every continuous $h\colon[0,T]\to \R$.} to a finite\footnote{If we were working on the unbounded time interval $[0,\infty)$ instead of $[0,T]$, we would have to ask, following \cite{Follmer1981}, that $\mu$ is Radon  (i.e., finite on compacts) and that $\mu_n\to \mu$ vaguely (i.e., $\int h \dd \mu_n\to \int h \dd \mu$ for every continuous  $h$ \emph{with compact support}).} measure $\mu$ such that the jumps of the (increasing, c{\`a}dl{\`a}g) function $[x]_t:=\mu([0,t])$ are given by $\Delta [x]_t=(\Delta x_t)^2$ for all $t\in [0,T]$. $\mathbb{Q}((\pi^n)_n)$ denotes the set of functions in $D([0,T];\R)$ having a quadratic variation along $(\pi^n)_n$. 
\end{definition}

For $x\in \mathbb{Q}((\pi^n)_n)$, we write $[x]^c$ and $[x]^d$ for the continuous and purely discontinuous parts of the c{\`a}dl{\`a}g function $[x]$ and note that by the above definition we have
\begin{equation*}
  \textstyle [x]^d_t = \sum_{0<s\leq t} (\Delta x_s)^2,\quad 0< t\leq T.
\end{equation*}
We now recall F{\"o}llmer's pathwise version of It{\^o}'s formula for paths in $\mathbb{Q}((\pi^n)_n)$. Here and throughout, $\int_0^t$ stands for $\int_{(0,t]}$ and increasing is understood as non-decreasing. 
\begin{theorem}[\cite{Follmer1981}]\label{thm:ito formula}
  Let $x \in \mathbb{Q}((\pi^n)_n)$ and $f \in C^2(\mathbb{R})$. Then, the pathwise It\^o formula
  \begin{align}\label{ito formula}
    f(x_t) -f(x_0)  = \int_0^t  f^\prime(x_{s-}) \dd x_s + \frac{1}{2} \int_0^t f^{\prime \prime} (x_s) \dd [x]^c_s 
    + J_t^f(x), \quad t\in [0,T], 
  \end{align} 
  holds with $J_t^f(x)$ as in~\eqref{eq:sum of jumps}, and with   
  \begin{equation}\label{eq:int for smooth integrands}
    \int_0^t f^\prime (x_{s-}) \dd x_s := \lim_{n \to \infty} \sum_{t_j \in \pi^n(t) } f^\prime (x_{t_j}) (x_{t_{j+1}} - x_{t_j}), \quad t \in [0,T],
  \end{equation}
  where the series in~\eqref{eq:sum of jumps} is absolutely convergent and the limit in~\eqref{eq:int for smooth integrands} exists.
\end{theorem}

We note that, to define $\int_0^t f^\prime(x_{s-}) \dd x_s $, F{\"o}llmer~\cite{Follmer1981} takes limits of sums of the form 
\[
  \sum_{\pi^n \ni t_j\leq t} g(x_{t_j})(x_{t_{j+1}}-x_{t_j}) , 
  \text{ whereas we consider }
  \sum_{t_j \in \pi^n}g(x_{t_j})(x_{t_{j+1}\wedge t}-x_{t_{j}\wedge t}).
\]
This however has no consequences, since the difference between these two sums is
$$ 
  g(x_{t_c(\pi^n,t)}) (x_{t_{c(\pi^n,t)+1}}-x_t),  \quad  \text{where } c(\pi,t):=\max\{j: \pi \ni t_j \leq  t\},
$$
which goes to zero as $|\pi^n|\to 0$ since $g$ is bounded on $[\inf_{t\in [0,T]} x_t , \sup_{t\in [0,T]} x_t]$, $x$ is c{\`a}dl{\`a}g and $t<t_{c(\pi,t)+1}\leq t+|\pi|$. In consequence, F{\"o}llmer's pathwise It{\^o} formula~\eqref{ito formula} holds also with our definition of $\int_0^t f^\prime(x_{s-}) \dd x_s $ and we shall exploit it in our proofs. Notice that analogously  
\begin{align*}
  \sum_{\pi^n \ni t_j\leq t} g(x_{t_j})(x_{t_{j+1}}-x_{t_j})^2  
  \quad \text{and} \quad  
  \sum_{t_j \in \pi^n}g(x_{t_j})(x_{t_{j+1}\wedge t}-x_{t_{j}\wedge t})^2
\end{align*} 
differ by
\[
  g(x_{t_c})((x_{t_{c+1}}-x_{t_c})^2 - (x_{t}-x_{t_c})^2 )  , \quad\text{with } c= c(\pi^n,t),
\]
which goes to zero as $|\pi^n| \to 0$.

\subsection{Local time via occupation measure}

In order to extend the It{\^o} formula for twice continuously differentiable functions~$f$ to twice \textit{weakly} differentiable functions~$f$, the notion of quadratic variation is not sufficient and the concept of local time is required. In probability theory there exist various classical approaches to define local times of stochastic processes. In the present deterministic setting, we first introduce a pathwise local time corresponding to the notation of local time as an occupation measure with respect to the quadratic variation. 

\medskip

The space of $q$-integrable (equivalence classes of) functions $g\colon \R\to \R$ is denoted by $L^q(\R):=L^q(\R;\R)$ with corresponding norm $\|\cdot \|_{L^q}$ for $q\in [1,\infty]$ and $W^{k,q}(\R) := W^{k,q}(\R;\R)$ stands for the Sobolev space of functions $g\colon \R\to \R$ which are $k$-times weakly differentiable in $L^q(\R)$, for $k\in \N$. Moreover, $L^q(K;\R)$ is the space of $q$-integrable functions $f\colon K\to \R$ for a Borel set $K\subset \R$ and we recall the left-continuous sign-function
\begin{equation*}
  \sign (x) := \begin{cases}
                 1&\text{if }x>0\\
                 -1&\text{if } x\leq 0
               \end{cases}.
\end{equation*}
We define, for $a,b \in \mathbb{R}$, 
\begin{align*}
  \llbracket a,b \rrparenthesis  := \begin{cases}
                                     [a,b) & \text{if } a \leq b\\
                                     [b,a) & \text{if } a > b
                                   \end{cases}
                                   \quad \text{with}\quad [a,a):=\emptyset.
\end{align*}

\begin{definition}
  Let $x \in \mathbb{Q}((\pi^n)_n)$. A Borel function $L_\cdot(x,\cdot) \colon [0,T] \times \mathbb{R} \to [0,\infty)$ is called the \emph{occupation local time} of $x$ if 
  \begin{equation}\label{eq:odf}
    \int_{-\infty}^{\infty} g(u) L_t(x,u) \dd u =  \int_0^t g(x_s) \dd [ x ]^c_s , \quad t \in [0,T],
  \end{equation}
  holds for any positive Borel function $g\colon\R \to [0,\infty)$.
\end{definition}
Naturally, this approach to local time is not new, see for example \cite{Bertoin87}. 
To extend It{\^o}'s formula to a Tanaka--Meyer formula, as, e.g., in \cite{Protter2004}, we will consider the quantity 
$$J_t (x,\cdot):=J_t^{f_u}(x),\quad \text{where }f_u:=|\cdot - u|/2.$$ We will, at times, drop $x$ from the notation, and simply write $L_t(u)$ and $J_t(u)$. It is straightforward to verify\footnote{Either checking separately the six cases where $u\leq x_{s-}\leq x_s$, $x_{s-}\leq u \leq x_s$ etc., or using the identity~\eqref{eq:J^fab=intf} with the function $f_u(\cdot):=|\cdot - u|/2$ and noting that $f_u'(\cdot)=\textup{sign}(\cdot - u)$.} that
\begin{align}\label{eq:Two expressions for J}
  |x_s-u| - |x_{s-} -u|- \textup{sign}(x_{s-}-u )\Delta x_s=  2 |x_s-u| \1_{  \llbracket x_{s-}, x_{s} \rrparenthesis  } ,
\end{align} 
which yields the useful compact expression 
\begin{align}\label{eq: def J_t(x,u)} 
  \textstyle  J_t (x, u)=\sum_{0<s\leq t} |x_s-u| \1_{  \llbracket x_{s-}, x_{s} \rrparenthesis  }(u) ,\quad u\in \R ,
\end{align} 
which readily implies that $J$ is a positive and increasing function. In particular, see Remark \ref{rem:LJcadlag} below, $L_t(\cdot)/2+ J_t (\cdot)\in L^p(\R)$ if and only if $L_t(\cdot), J_t (\cdot)\in L^p(\R)$. Notice that $x$ is bounded, since it is c{\`a}dl{\`a}g, and $L_t(u)$ and $J_t ( u)$ equal $0$ if $u$ does not belong to the compact set $[\inf_{s\in [0,T]} x_s, \sup_{s\in [0,T]} x_s ]$. 

\begin{definition}\label{def:L^pLT}
  We let $\mathbb{L}_p((\pi^n)_n)$ denote the set of all paths $x\in \mathbb{Q}((\pi^n)_n)$ having an occupation local time $L$ and such that $K_t(x,\cdot):=L_t(x,\cdot)/2+ J_t (x,\cdot)\in L^p(\R)$ for all $t \in [0,T]$.
\end{definition} 

There is no common agreement in the related literature in probability theory as to whether $L$ or $L/2$ is to be called local time, cf. \cite[Remark~6.4]{Karatzas1988}; here we decided to follow the convention made in the standard textbook~\cite{Protter2004}. A classical approach to extend It\^o's formula and, in particular, the ``stochastic'' integral $\int_0^t f^\prime(x_{s-}) \dd x_s$ to twice weakly differentiable functions~$f$, is to approximate the function $f$ by smooth functions, cf. \cite[Theorem~3.6.22]{Karatzas1988} for the case of Brownian motion. For this purpose we consider a ``mollifier'' $\rho$, i.e., a positive function $\rho \in C^{\infty}(\mathbb{R})$ and such that $\int_{-\infty}^\infty \rho(u) \dd u= 1$, and set $\rho_n(u):= n \rho (nu)$ for $n \in \mathbb{N}$. Given a function $f \in W^{2,q}(\R)$ we approximate it via the convolution $f_n := \rho_n * f$. In this way, $f_n\in C^2(\R)$, $f_n\to f$ in $W^{2,q}(\R)$ if $q<\infty$ (if $q=\infty$ this is true if one assumes $f''$ is continuous) and, in particular, $\lim_{n \to \infty} f_n (x) = f(x)$ for $x \in \mathbb{R}$.

\begin{proposition}\label{prop:ito formula via occupation measure}
  Let $x \in \mathbb{L}_p((\pi^n)_n)$ and $f\in W^{2,q}(\R)$ with $1/p+1/q\geq 1$ and $q\in [1,\infty)$. Then, the series~\eqref{eq:sum of jumps} defining $J_t^f(x)$ is absolutely convergent, $\int_0^t f_n^\prime (x_{s-}) \dd x_s$ defined by \eqref{eq:int for smooth integrands} converges to the finite limit 
  \begin{equation}\label{eq:intfdx:=lim_n}
    \int_0^t f^\prime (x_{s-}) \dd x_s := \lim_{n \to \infty} \int_0^t f_n^\prime (x_{s-}) \dd x_s, \quad t\in [0,T],
  \end{equation}
  which does not depend on the choice of $\rho$, and the pathwise Tanaka--Meyer formula
  \begin{equation}\label{eq:Tanaka formula}
    f(x_t)-f(x_0)  
    =\int_{0}^t f^\prime(x_{s-})\dd x_s+ \frac{1}{2} \int_{\R} L_t(x,u)  f^{''}(\mathrm{d}u) +J_t^f(x),\quad t\in [0,T], 
  \end{equation}
  holds with such definition of $\int_0^t f^\prime (x_{s-}) \dd x_s$.\\ The statements hold for $q=\infty$ if $f''$ is continuous.
\end{proposition}

Because of Proposition~\ref{prop:ito formula via occupation measure}, it is of interest to ask under which assumptions one can get that $L_t(x,\cdot)$ and $J_t(x,\cdot)$ are in $L^p(\R)$. First, remark that, since both quantities are equal to $0$ outside a compact, the $p$-integrability requirement in Definition~\ref{def:L^pLT} is a local one. Then, notice that if $x\in \mathbb{Q}((\pi^n)_n)$ has an occupation local time then $L_t,J_t \in L^1(\R)$ (i.e., $x\in \mathbb{L}_1((\pi^n)_n)$), since 
\begin{align*}
  \int_{\R} L_t(x,u) \dd u =  [ x ]^c_t <\infty , \qquad \int_{\R} J_t (x, u) \dd u = \frac{1}{2} [ x ]^d_t <\infty.
\end{align*} 

\begin{remark}\label{Jfinite} 
  If $p\in [1,\infty)$ and $C_p:=1/(p+1)^{1/p}$ then 
  $$
    \|J_t (x, \cdot)\|_{L^p} \leq C_p \sum_{0<s\leq t} |\Delta x_s|^{1+\frac{1}{p}} .
  $$
  This can be seen as a consequence of Minkowski's integral inequality and of the identity
  \begin{align}\label{eq:int |b-u|^p}
    \int_{  \llbracket a,b \rrparenthesis } |b-u|^p \dd u=\frac{|b-a|^{p+1}}{p+1} .
  \end{align} 
  A similar bound for $L$ can be given under the stronger assumption  $x \in \mathbb{L}_p^W((\pi^n)_n)$, see Definition~\ref{def:Wuermli local time} and equation~\eqref{eq:L^p bound for L} in the next subsection. Alternatively, if $x \in \mathbb{L}_1((\pi^n)_n)$, then $p$-summability for $L$, for $p\in (1,\infty)$, is equivalent to: 
  \begin{equation*}
    \|L_t(x,\cdot)\|_{L^p}=\sup \Big\{\int_0^t g(x_s) \dd [ x ]^c_s  :  \|g\|_{L^q}\leq 1 \Big\}<\infty.
  \end{equation*}
\end{remark} 

Notice that an occupation local time~$L$ is only unique up to equality a.e.\footnote{Here, and elsewhere unless otherwise specified, a.e. $u$ is with respect to the Lebesgue measure.} $u$ for each $t$; in particular, $L$ could be thought of as an equivalence class, and one is then led to look for good representatives. In particular, it is often of interest to have a version $L$ which is c{\`a}dl{\`a}g in $t$. This can be ensured along the same lines as standard results on c{\`a}dl{\`a}g version of supermartingales since $L_s\leq L_t$ a.e.\ for any $0\leq s\leq t$, $L_t\in L^1(\R)$ for all $t$ and $t \mapsto \int_{\R} L_t(u) \dd u=[x]^c_t$ is continuous. Similarly, existence of a c{\`a}dl{\`a}g version for $J$ follows from the fact that $J_T(u)<\infty$ for a.e. $u$, that $x$ is c{\`a}dl{\`a}g and that $J_t(x,\cdot)$, see \eqref{eq: def J_t(x,u)}, is defined using jumps of $x$ up to and including time $t$. 

\begin{remark}\label{rem:LJcadlag} 
  If $x$ has an occupation local time $L$, then one can choose for each $t\in [0,T]$ a version $\tilde{L}_t(\cdot)$ of $L_t$ such that $\tilde{L}_\cdot(u)$ is positive, finite, c{\`a}dl{\`a}g and increasing for each $u\in \R$. Moreover, $J_\cdot(u)$ is positive, finite and c{\`a}dl{\`a}g increasing for a.e. $u$. In particular, it follows that $L_t(\cdot), J_t (\cdot)\in L^p(\R)$ holds for every $t\in [0,T]$ if and only if $L_T(\cdot), J_T (\cdot)\in L^p(\R)$.
\end{remark} 

It can also be useful to have right-continuity of $L,J$ in the variable $u$. For $J$ here is a simple criterion; for $L$, it has to be assumed: cf. Remark~\ref{PRemark:F} below. 

\begin{remark}\label{Jucadlag}
  Notice that 
  \begin{equation*}
    \TTV{J_t (x, \cdot)}{\R}{}
    := \sup \bigg \{ \sum_{i=0}^{N-1} |J_t (x, u_{i+1})-J_t (x, u_{i})|\,:\, (u_i)_{i=0}^N\subset \R,\, N\in\N  \bigg\} \leq \sum_{0<s\leq t }  |\Delta x_s| ,
  \end{equation*}
  and so if $\sum_{0<s\leq t } |\Delta x_s|<\infty$ for all $t$, then $J_t (x, \cdot)$ is c{\`a}dl{\`a}g and of finite variation for all $t\in [0,T]$.
\end{remark} 

As an application of having a version $\tilde{L}$ of $L$ which is c{\`a}dl{\`a}g in $t$, notice that the occupation time formula~\eqref{eq:odf} then extends to all positive Borel $h=h(s,u)$ as follows
\begin{align*}
  \int_{-\infty}^{\infty} \left(\int_0^t  h(s,u) \dd \tilde{L}_s(x,u)\right) \dd u = \int_0^t h(s,x_s) \dd [ x ]^c_s ,\quad t\in [0,T].
\end{align*} 
Moreover, since $J$ is c{\`a}dl{\`a}g in $t$ it also satisfies a restricted occupation time formula: if $h=h(s,u)$ is a positive Borel function such that $h(s,u)=h(s,x_s)$ for a.e.\ $u\in \llbracket x_{s-}, x_{s} \rrparenthesis$, then Fubini's theorem gives that 
\begin{align*}
  \int_{-\infty}^{\infty} \left(\int_0^t  h(s,u) \dd J_s(x,u)\right) \dd u =  \frac{1}{2} \int_0^t h(s,x_s) \dd [ x ]^d_s ,
\end{align*} 
and this observation seems to be new. 

\medskip

To facilitate the proof of Proposition~\ref{prop:ito formula via occupation measure}, as well as for later use, let us recall some well known facts. A function $g\colon\R\to \R$  is convex iff its second distributional derivative $g''$  is a positive Radon measure. Thus $f\colon\R\to \R$ equals to the difference of two convex functions iff $f''$ is a signed Radon measure. We may then write $f=g-h$ with $g,h$ convex and $|f''|=g'' + h''$ being the measure associated with the total variation of $f$, $ \TTV{f'(\cdot)}{[0,t]}{}=|f''|([0,t])$. Given such $f$, $f'$ denotes the left-derivative of $f$, which is left-continuous and of locally bounded variation and satisfies $f(b)-f(a) = \int_a^b f^\prime(y) \dd y$ for all $a,b\in \R$. Thus for $b\ge a$ we get that
\begin{equation*}
  f(b) - f(a) - f^\prime(a)(b-a) = \int_a^b (f^\prime(u) - f^\prime(a)) \dd u =  \int_{[a,b)} (b-u)  \,f''(\mathrm{d} u),
\end{equation*}
where we used integration by parts. For $b<a$, we get instead 
\begin{equation*}
  f(b) - f(a) - f^\prime(a)(b-a) = \int_{[b,a)} (u-b)\, f''(\mathrm{d} u),
\end{equation*}
so 
we obtain the identity
\begin{align}\label{eq:J^fab=intf}
  J^f(a,b):=  f(a)- f(b)-f^\prime (b)(a-b)   
  = \int_{  \llbracket a,b \rrparenthesis } |b-u| \, f''(\mathrm{d} u), \qquad a,b \in \R ,
\end{align} 
which can often be used in proofs in lieu of the following representation 
\begin{equation}\label{eq:representation of f}
  f(x)= a x+ b +  (|\cdot | * f^{\prime\prime}) (x), \quad x\in \R ,
\end{equation}
(which holds for some  $a,b\in\R$), which is often used in the literature. Representation~\eqref{eq:representation of f} holds whenever $\int_\R |a-u|\, f''(\d u)<\infty$ for all $a$ (in particular if $f''$ has compact support), and is proved after Proposition~3.2 in \cite[Appendix~3]{Revuz1999}.

A version of the following statement appears without proof during the course of the proof of \cite[Chapter~4, Theorem~70]{Protter2004}.

\begin{lemma}\label{JfJut}
  If $f\colon\R\to \R$ is a convex function then the series~\eqref{eq:sum of jumps} defining $J_t^f(x)$ consists only of positive terms. If $f$ equals to the difference of two convex functions and 
  \begin{align*}
    \int_{\R}  J_t(x,u)\,  |f''|(\mathrm{d} u)<\infty,
  \end{align*} 
  then the series~\eqref{eq:sum of jumps} is absolutely convergent. In both cases, the series~\eqref{eq:sum of jumps} defining $J_t^f(x)$ is well defined\footnote{See \eqref{Sum AC}.} and satisfies
  \begin{align}\label{J_tint}
    J_t^f(x)=  \int_{\R}  J_t(x,u) \, f''(\mathrm{d} u),\quad t\in [0,T].
  \end{align} 
\end{lemma} 

\begin{proof}
  From~\eqref{eq:J^fab=intf} we get 
  \begin{align}\label{eq:J^feqint}
    J^f(x_s,x_{s-}) =\int_{\R}  |x_s-u| \1_{  \llbracket x_{s-}, x_{s} \rrparenthesis  }(u)\, f''(\mathrm{d} u) .
  \end{align} 
  If $f$ is convex the series~\eqref{eq:sum of jumps} defining $J_t^f(x)$ consists only of positive terms, and the thesis follows from \eqref{eq:J^feqint}, summing  over $s\leq t$  and applying Fubini's theorem. If instead $f=g-h$ with $g,h$ convex then $|f''|=g''+h''$ and, by assumption, $\int_{\R} J_t(x,u) \, |f''|(\mathrm{d} u)<\infty$. \eqref{J_tint} follows again from Fubini's theorem. The absolute convergence of the series~\eqref{eq:sum of jumps} follows writing 
  \begin{equation*}
    |\Delta f(x_{s})-f^\prime (x_{s-} )\Delta x_s | \leq  (\Delta g(x_{s})-g^\prime (x_{s-} )\Delta x_s) + (\Delta h(x_{s})-h^\prime (x_{s-} )\Delta x_s) , 
  \end{equation*}
  summing the latter over $s\leq t$ and applying~\eqref{J_tint} to $g$ and $h$.
\end{proof} 

\begin{remark}\label{rem: J^f finite}
  It follows from Lemma~\ref{JfJut} and H{\"o}lder's inequality that, if $J_t(x,\cdot)\in L^p(\R)$ and $f''(\d u)=f''(u) \dd u$ with $f''\in L^q(\R)$, where $p,q\geq 1$ are conjugate exponents, that is satisfy $1/p+1/q = 1$, then the series~\eqref{eq:sum of jumps} defining $J_t^f(x)$ is absolutely convergent. Moreover, if $J_t(x,\cdot)$ is  bounded\footnote{This happens for example if $\sum_{s\leq t} |\Delta x_s|<\infty$, by Remark~\ref{Jucadlag}.} then the series~\eqref{eq:sum of jumps} is absolutely convergent for every $f$ which is a difference of convex functions: indeed, $J_t(x,\cdot)=0$ outside a compact, and $|f''|(C)<\infty$ for every compact $C\subseteq \R$.
\end{remark} 

An alternative, possibly more intuitive but also more cumbersome, way of getting~\eqref{eq:J^feqint} is to define
$$
  g(\cdot):=|x_s- \cdot | \1_{  \llbracket x_{s-}, x_{s} \rrparenthesis  }(\cdot ) ,
$$ 
which is in $ L^1(\R)$, equals zero outside a compact, and has distributional derivatives
$$
  \D g=(\Delta x_s) \delta_{x_{s-}}-\1_{[x_{s-},\infty)}+\1_{[x_s,\infty)}, \quad \D^2 g=(\Delta x_s) \D\delta_{x_{s-}} -\delta_{x_{s-}}+\delta_{x_{s}} .
$$
Then equation~\eqref{eq:J^feqint} is simply\footnote{This equality holds a priori only when $f$ is $C^{\infty}(\R)$ (by definition of distributional derivatives). However, with some work it follows that it holds for any $f$ which equals the difference of convex functions: indeed, since $g$ is c{\`a}dl{\`a}g, convolving against a mollifier with support in $[0,\infty)$ shows that there exist $f_\epsilon\in C^\infty(\R)$ such that $f_\epsilon\to f$ uniformly on compacts and $\int g(u) (\D^2 f_\epsilon)(u) \dd u \to \int  g(u) \dd (\D^2 f)$, as shown in \cite[Proof of Theorem~5.2]{Davis2018}.} the identity $ \int_{\R}  f(u)\, (\D^2 g)(\mathrm{d} u)=\int_{\R} g(u) \,(\D^2 f)(\mathrm{d} u)$.
  
\begin{proof}[Proof of Proposition~\ref{prop:ito formula via occupation measure}] 
  The series~\eqref{eq:sum of jumps} defining $J_t^f(x)$ is absolutely convergent by Remark~\ref{rem: J^f finite}. If $h\in C^2(\R)$, from F{\"o}llmer's pathwise It\^o formula~\eqref{ito formula}, the definition of occupation local time $L$ and of $K:=L/2+J$, and Lemma~\ref{JfJut}, it follows that 
  \begin{equation}\label{eq:TanakaformulaInt}  
    \int_{0}^t h^\prime(x_{s-})\dd x_s =
    h(x_t)-h(x_0)  -  \int_{\R}  K_t(x,u) \, h''(\mathrm{d} u),\quad t\in [0,T],
  \end{equation}
  holds with $\int_0^t h^\prime (x_{s-}) \dd x_s $ defined via~\eqref{eq:int for smooth integrands}. Applying~\eqref{eq:TanakaformulaInt} to $h=f_n= \rho_n * f  \in C^2(\R)$ and taking limit as $n\to \infty$, the right-hand side converges to 
  $$
    f(x_t)-f(x_0)  -  \int_{\R}  K_t(x,u)  f''(u) \,\mathrm{d} u 
  $$
  because $K_t(x,\cdot)\in L^p(\R)$ and $f_n\to f$ in $W^{2,q}(\R)$ (so  $f_n\to f$ pointwise and $f_n''\to f''$ in $L^q(\R)$). It follows that the LHS converges as well.
\end{proof}

\begin{remark}
  It follows from~\eqref{J_tint} that, whenever Tanaka--Meyer's formula holds, it can be written as 
  \begin{equation} \label{eq:Tanaka formulaNoJf}
    f(x_t)-f(x_0) = \int_{0}^t f^\prime(x_{s-})\dd x_s +\int_{\R}  K_t(x,u)\,  f''(\mathrm{d} u) ,\quad t\in [0,T]  ,
  \end{equation}
  where we recall that $K_t(u) :=L_t(u)/2 + J_t(u)$. While uncommon, writing \eqref{eq:Tanaka formulaNoJf} seems rather elegant and simpler than \eqref{eq:Tanaka formula}.
\end{remark}

\begin{remark}\label{PRemark:C}
  One can recover a continuous in time, for a.e. level $u$, version $\tilde{L}$ of the occupation time $L$ from knowing just a jointly measurable function $K_t(u)$ such that $K_\cdot(u)$ is c{\`a}dl{\`a}g increasing for a.e. $u$, $K_0=0$, $K_T\in L^1(\R)$, and \eqref{eq:Tanaka formulaNoJf} holds for all $f\in C^2$ with $\int_0^t f^\prime (x_{s-}) \dd x_s $ defined via \eqref{eq:int for smooth integrands}. Indeed, {$\tilde{L}_\cdot(u)$ (resp. $J_\cdot(u)$) is the continuous (resp. purely discontinuous) part of the increasing c{\`a}dl{\`a}g function $K_\cdot(u)$}. To show this, consider that for $f\in C^2(\R)$ F{\"o}llmer's formula~\eqref{ito formula}, \eqref{eq:TanakaformulaInt} and Lemma~\ref{JfJut} give that 
  \begin{equation*}
    K^f_t
    :=\int_{\R}  K_t^c(u)  f'' (u) \,\mathrm{d} u +\int_\R K_t^d(u)  f''(u)\, \mathrm{d} u 
    = \frac{1}{2} \int_0^t f^{\prime \prime} (x_s) \dd [x]^c_s +\int_{\R}  J_t(u) f''(u)\, \mathrm{d}  u ,
  \end{equation*}
  where $K^c$ (resp. $K^d$) denotes the continuous (resp. purely discontinuous) part of $K_\cdot(u)$. In each of the two above representations of the c{\`a}dl{\`a}g increasing function $K^f_t$ the first term is continuous and the second purely discontinuous, so by uniqueness of such decomposition 
  $$ 
    \int_{\R}  K_t^c(u)  g( u)\, \mathrm{d} u= \frac{1}{2} \int_0^t g(x_s) \dd [x]^c_s, \qquad  \int_\R K_t^d(u)  g(u)\, \mathrm{d}  u 
    =  \int_{\R}  J_t(u) g(u)\, \mathrm{d} u 
  $$
  holds for any $g$ of the form $f''$, i.e., for any continuous $g$; but then it also automatically holds for any Borel $g$, so $2K^c$ is an occupation local time of $x$ and $J_t=K^d_t$ a.e. $u$ for each $t$; since $J_t$ and $K^d_t$ are c{\`a}dl{\`a}g in $t$, $J_t=K^d_t$ a.e. $u$ for all $t$.
\end{remark}
 
\begin{remark}
  For continuous paths~$x$ the above approximation argument can be used to obtain space-time Tanaka--Meyer formulae without relying on the representation~\eqref{eq:representation of f}, see \cite{Feng2006}. Although elaborated in a probabilistic framework, the proofs in \cite{Feng2006} are (primarily) of pathwise nature.
\end{remark}

\begin{remark}\label{rem:pathwise Ito formula}
  The definition of occupation local times and the generalization of It{\^o}'s formula to only twice weakly differentiable functions in Proposition~\ref{prop:ito formula via occupation measure} is based on F{\"o}llmer's notion of quadratic variation and his pathwise It\^o formula (Theorem~\ref{thm:ito formula}). However, the F{\"o}llmer--It{\^o} formula is by no means the only pathwise It{\^o}-type formula, which can be extended to an Tanaka--Meyer formula in the spirit of in Proposition~\ref{prop:ito formula via occupation measure}. For example, one could also start from the pathwise It{\^o} formula based on c{\`a}dl{\`a}g rough paths (\cite[Theorem~2.12]{Friz2018}) or the one based on truncated variation (\cite[Theorem~4.1]{Lochowski2019}) and proceed in an analogous manner as done in the present subsection.
\end{remark}

\begin{remark}\label{PRemark:F}
  If $x$ has an occupation local time $L$, then one can give explicit formulae for $L$. Indeed, since $L_t(\cdot)\in L^1(\R)$, taking $\lim_{\epsilon \downarrow 0}$ of \eqref{eq:odf} applied to $g:=\1_{[u-\epsilon, u+\epsilon]}$ gives that 
  $$ 
    L_t(x,u)= \lim_{\epsilon \downarrow 0} \frac{1}{2\epsilon} \int_0^t \1_{[u-\epsilon, u+\epsilon]}(x_s) \dd [ x ]^c_s , \quad\text{for a.e. } u,
  $$
  meaning that the limit on the right-hand side exists for a.e.\ $u \in \R$  and is a version of $L_t(\cdot)$. 
  Analogously, if we can apply Tanaka--Meyer's formula to the convex function $|\cdot - u|$ we get the following expression for $L$:
  \begin{equation}\label{eq:def L using | |}
    L_t(x,u)=|x_t - u|-|x_0 - u| -\int_{0}^t \sign(x_{s-}-u)\dd x_s- 2 J_t(x,u),\quad t\in [0,T]. 
  \end{equation}
  It is thus desirable to establish if (a version of) Proposition~\ref{prop:ito formula via occupation measure} holds in the case where $f\colon \R\to \R$ equals to the difference of two convex functions. This is the case under the additional assumptions that the mollifier~$\rho$ has compact support in $[0,\infty)$, that $J_t(u)$ is c{\`a}dl{\`a}g in $u$ for all $t$ (see Remark~\ref{Jucadlag}), and that there exists a version $\tilde{L}_t$ of the pathwise local time $L_t$ which is c{\`a}dl{\`a}g in $u$ for all $t$ (in particular, unlike in the stochastic setting, one cannot use \eqref{eq:def L using | |} to prove that $L$ has a version which is c{\`a}dl{\`a}g in $u$ for all $t$ without running into circular arguments). Indeed, under these assumptions the proof of \cite[Theorem~5.2]{Davis2018} shows that $\int_{\R} g(u)\, f_n''(\d u) \to \int_{\R} g(u) \,f''(\d u) $ for any c{\`a}dl{\`a}g $g$, and if we apply this to $g=K_t$ the rest of the proof of Proposition~\ref{prop:ito formula via occupation measure} goes through.
\end{remark}

\begin{remark}
  As in\footnote{One can apply the proof found in \cite[Chapter~4, Theorem~69]{Protter2004}, which simplifies somewhat as we do not need to deal with the dependence on $\omega$.} the stochastic setting, if we can\footnote{See Remark~\ref{PRemark:F}.} apply Tanaka--Meyer's formula to the convex function $f(x)=(x-u)^+$, we find that the measure  $\dd L_\cdot(u)$ is supported by the set $\{s\in (0,t]:x_{s}=x_{s-}=u\}$, and correspondingly, the measure $\dd J_\cdot(u)$ is carried by the set 
  $$
    \{s\in (0,t]:u\in ( x_{s-}, x_{s} ) \text{ or } u\in ( x_{s}, x_{s-} ] \}  
  $$ 
  of times at which $x$ jumps across\footnote{More precisely, if the jump is downward, then $x$ is allowed to jump from $x_{s-}=u$.} level $u$.
\end{remark}

\subsection{Local time via discretization}

An alternative approach to achieve a pathwise Tanaka--Meyer formula goes back to W{\"u}rlmi~\cite{Wuermli1980} and is based on a discrete version of the Tanaka--Meyer formula. For continuous paths~$x$ this approach is well-understood and led to several extensions, see \cite{Perkowski2015,Davis2018,Cont2019}. One feature of this discretization argument is that the ``stochastic'' integral~$\int_0^t f^\prime(x_{s-}) \dd x_s$ is still given as a limit of left-point Riemann sums, see also \cite{Davis2014}. In the present subsection we generalize W{\"u}rlmi's approach to the case of c{\`a}dl{\`a}g paths~$x$. Given a partition $\pi=(t_j)_{j=0}^n$ of $[0,T]$, we define the discrete level crossing time of $x$ at $u$ (along $\pi$) as the function
\begin{align}\label{eq:Kpi}
  \textstyle K^{\pi}_t(x,u):=\sum_{t_j\in \pi} |x_{t_{j+1}\wedge t} -u|\1_{ \llbracket x_{{t_j}\wedge t}, x_{{t_{j+1}\wedge t}} \rrparenthesis}(u), \quad t\in [0,T].
\end{align} 
Then, applying \eqref{eq:J^fab=intf} to $a=x_{t_i \wedge t}, b=x_{t_{i+1} \wedge t}$ and summing over $i$, we obtain the discrete version of Tanaka--Meyer formula 
\begin{align}\label{eq:discrete tanaka-meyer formula}
  f(x_t) - f(x_0)  -\sum_{t_i \in \pi} f^\prime (x_{t_i}) (x_{t_{i+1}\wedge t} - x_{t_i \wedge t}) 
  = \int_{\R} K^{\pi}_t(u)\, f^{\prime\prime}(\d u) .
\end{align}
Taking limits along a sequence of partitions $(\pi^n)_n$, with $|\pi^n|\to 0$, we obtain the following definition of $L^p$-level crossing time. We note that it extends the previous works for continuous paths, e.g., \cite[Definition~B.3]{Davis2014}. We also note that using the same notation $K_t$ as before will be justified by Proposition~\ref{prop:meyer-ito formula}.

\begin{definition}\label{def:Wuermli local time}
  Let $x \in D([0,T];\mathbb{R})$ and let $(\pi^n)_n$ be a sequence of partitions such that $|\pi^n|\to 0$. A function $K \colon [0,T] \times \mathbb{R} \to \mathbb{R}$ is called the $L^p$-\textup{level crossing time} of $x$ (along $(\pi^n)_n$) if $ K^{\pi^n}_t$ converges weakly in $L^p(\R)$ to $K_t$ for each $t\in [0,T]$ as $n \to \infty$, and $t \mapsto \int_\R K_t(u) \dd u$ is right-continuous. The set $\mathbb{L}^{W}_p((\pi^n)_n)$ denotes all paths $x\in D([0,T];\mathbb{R})$ having an $L^p$-level crossing time along~$(\pi^n)_n$. 
\end{definition}

\begin{lemma}
  The level crossing time $K$ in Definition~\ref{def:Wuermli local time} is increasing in $t\in [0,T]$, i.e., $K_s(\cdot)\leq K_t(\cdot)$ a.e. for each $s\leq t$. 
\end{lemma} 

\begin{proof}
  Given $\pi=(t_j)_j$, let $m(\pi,s)$ be the value of $j$ such that $t_j<s\leq t_{j+1}$, and write 
  $$
    K^\pi_s=\sum_{j < m(\pi,s)}  a_j(u) + |x_{s} -u|\1_{\llbracket x_{t_{m(\pi,s)}},x_s \rrparenthesis}(u)  , \text{ where } a_j(u):=|x_{t_{j+1}} -u|\1_{ \llbracket x_{{t_j}}, x_{{t_{j+1}}} \rrparenthesis}(u).
  $$
  If $s< t$, analogously write 
  $$
    K^\pi_t-\sum_{j < m(\pi,s)}  a_j(u)- a_{m(\pi,s)}(u)= \sum_{m(\pi,s)<j < m(\pi,t)}  a_j(u) +|x_{t} -u|\1_{ \llbracket x_{t_{m(\pi,t)}},x_t \rrparenthesis}(u) =:R_{s,t}^\pi .
  $$
  Thus
  $$
    K^\pi_t- K^\pi_s-R_{s,t}^\pi=a_{m(\pi,s)}(u)-|x_{s} -u|\1_{ \llbracket x_{t_{m(\pi,s)}},x_s \rrparenthesis}(u)  =:S_s(\pi,u) ,
  $$
  and since $R_{s,t}^\pi\geq 0$ the thesis follows once we prove that $S_s(\pi^n,u) \to 0$ for every $u$ when $|\pi^n|\to 0$. This holds since if $m(n):=m(\pi^n,s)$ then $t_{m(n)}$ and $ t_{m(n)+1}$ converge to $s$, and $t_{m(n)}< s \leq t_{m(n)+1}$, so 
  $$
    a_{m(n)}(u) \quad \text{ and } \quad |x_{s} -u|\1_{\llbracket x_{t_{m(n)}},x_s\rrparenthesis}(u) 
  $$
  both converge to $|x_{s} -u|\1_{\llbracket x_{s-},x_s \rrparenthesis}(u) $ as $n\to \infty$, since $x$ is c{\`a}dl{\`a}g.
\end{proof} 

Notice that $K_t$ is only defined as an equivalence class. Using the same arguments as in the discussion preceding Remark~\ref{rem:LJcadlag}, for each $t$ we can take the version of $K_t$ such that the resulting process is c{\`a}dl{\`a}g increasing in $t$ for each $u$. From now on, we will always work with such a version and we let $K^c$ (resp. $K^d$) denote the continuous (resp. purely discontinuous) part of the increasing c{\`a}dl{\`a}g function $K_\cdot(u)$. 

\begin{proposition}\label{prop:meyer-ito formula}
  Suppose that $x\in\mathbb{L}^{W}_p((\pi^n)_n)$ for $p,q\in [1,\infty]$, with $1/p+1/q=1$, and let $K$ be the $L^p$-\textup{level crossing time} of $x$ along $(\pi^n)_n$.
  If $f\in W^{2,q}(\R)$, then the following limit exists (and is finite)
  \begin{equation}\label{eq:stochastic integral as left-point Riemann sums}
    \int_0^t f^\prime(x_{s-}) \dd x_s := \lim_{n \to \infty}\sum_{t_i \in \pi^n} f^\prime (x_{t_i}) (x_{t_{i+1} \wedge t} - x_{t_i\wedge t}),\quad t\in [0,T],
  \end{equation}
  and the pathwise Tanaka--Meyer formula~\eqref{eq:Tanaka formulaNoJf} holds with such definition of $\int_0^t f^\prime (x_{s-}) \dd x_s$ and $K$. Moreover, $2K^c$ is the occupation local time of $x$ and $J_t(u)$ in \eqref{eq: def J_t(x,u)} satisfies $J_t(u)=K^d_t(u)$ for a.e. $u$ and for all $t\leq T$. In particular, also the pathwise Tanaka--Meyer formula~\eqref{eq:Tanaka formula} holds (with $L=2K^c$), the two definitions \eqref{eq:intfdx:=lim_n} and \eqref{eq:stochastic integral as left-point Riemann sums} of $ \int_0^t f^\prime (x_{s-}) \dd x_s $ coincide, and $\mathbb{L}^{W}_p((\pi^n)_n) \subseteq \mathbb{L}^{p}((\pi^n)_n)$.
\end{proposition}

\begin{proof}[Proof of Proposition~\ref{prop:meyer-ito formula}.]
  Taking the limit as $n$ goes to $\infty$ of the discrete Tanaka--Meyer formulae~\eqref{eq:discrete tanaka-meyer formula} applied to $\pi^n$, the RHS converges and hence also does the LHS. The pathwise Tanaka--Meyer formula \eqref{eq:Tanaka formulaNoJf} thus holds if using the definition \eqref{eq:stochastic integral as left-point Riemann sums}. Now, from Remark~\ref{PRemark:C} it follows that $2K^c$ satisfies the occupation time formula and the remaining statements readily follow. 
\end{proof}

\begin{remark}
  Following the seminal paper~\cite{Follmer1981}, we consider the ``stochastic'' integral as limit of left-point Riemann sums~\eqref{eq:stochastic integral as left-point Riemann sums} and not as limit of
  \begin{equation*}
    \sum_{t_i \in \pi^n} f^\prime (x_{t_i-}) (x_{t_{i+1} \wedge t} - x_{t_i\wedge t}), \quad t\in [0,T].
  \end{equation*}
  In a probabilistic setting, where $x$ is assumed to be a semimartingale, these limits coincide with the classical It{\^o} integral almost surely (see \cite[Chapter~II.5, Theorem~21]{Protter2004}) and so they are equal. In the present pathwise setting however, they could be different.
\end{remark}

\begin{remark}
  Applying Minkowski's integral inequality and using the identity~\eqref{eq:int |b-u|^p}, we obtain that if $p\in [1,\infty)$ and $C_p:=1/(p+1)^{1/p}$, then 
  \begin{align*}
    \|K_t^\pi \|_{L^p} \leq C_p \sum_{t_i\in \pi} |x_{t_{i+1}\wedge t} - x_{t_i \wedge t}|^{1+\frac{1}{p}} .
  \end{align*} 
  In particular, if $x \in \mathbb{L}^{W}_p((\pi^n)_n)$, then the occupation local time $L$ equals $2K^c$ and so satisfies 
  \begin{align}\label{eq:L^p bound for L}
    \|L_t \|_{L^p} \leq 2 \|K_t \|_{L^p} \leq 2 C_p \liminf_n \sum_{t_i\in \pi^n} |x_{t_{i+1}\wedge t} - x_{t_i \wedge t}|^{1+\frac{1}{p}}  \quad \text{ for every $p\in [1,\infty)$. }
  \end{align} 
\end{remark}

\begin{remark}
  Given the definition of $J_t(u)$, it seems natural that, if $x\in\mathbb{L}^{W}_p((\pi^n)_n)$ and
  \begin{equation*}
    J^{\pi}_t(u):= \sum_{t_i \in \pi(t)} \1_{\llparenthesis x_{t_i-}, x_{t_{i}} \rrbracket} (u) \vert x_{t_{i}} - u \vert ,\quad u\in \R ,\,t\in[0,T], 
  \end{equation*}
  then $J^{\pi^n}_t$ also converges weakly in $L^p(\R)$. If we assume this and  denote by $L^d_t$ the limit, if $(\pi^n)_n$ are refining and $\cup_n \pi^n \supseteq \{s\in [0,T]\,:\,\Delta x_s \neq 0\}$, then $L^d_t=J_t=K^d_t$ a.e. In particular, $K^{\pi^n}_t- J^{\pi^n}_t$ converges weakly in $L^p(\R)$ to $K^c_t$. Indeed, if $f'' \in L^q(\R)$, \eqref{eq:J^fab=intf} gives
  \begin{equation}\label{eq:JnfintJnu}
    J_t^{f,\pi^n}:=\sum_{t_i \in \pi^n(t)} f(x_{t_i}) - f(x_{t_{i} -})- f^\prime (x_{t_i-}) (x_{t_i} - x_{t_i-} ) 
    = \int_{\R} J_t^{\pi^n}(u) f^{\prime\prime}(u)\dd u ,
  \end{equation}
  so our assumptions and Lemma~\ref{JfJut} imply that the series \eqref{eq:sum of jumps} defining $J^f_t$ is absolutely convergent. Using the dominated convergence theorem we conclude that $ J_t^{f,\pi^n}\to J^f_t$, so taking $n\to \infty$ in \eqref{eq:JnfintJnu} we get  
  \begin{equation*}
    J_t^f=\int_{\R} L^d_t(u) f^{\prime\prime}(u)\dd u ,
  \end{equation*}
  so by Lemma~\ref{JfJut} $L^d_t=J_t$ a.e.
\end{remark}

\subsection{Local time via normalized numbers of interval crossings}

In Proposition~\ref{prop:ito formula via occupation measure} above we approximated $f$ with regular functions $f_n$ for which the ``stochastic'' integral $\int_0^t f_n^\prime(x_{s-}) \dd x_s$ was defined via Theorem~\ref{thm:ito formula}. An alternative regularisation idea would be to approximate the path $x$ by sufficiently regular functions~$(x^n)$, ensuring that the ``stochastic'' integral $\int_0^t f^\prime(x_{s-}^n) \dd x_s$ is well-defined for each $x^n$ (via integration by parts). We pursue this approach now using for $x^n$ the solutions to the so-called double Skorokhod problem. This choice of approximations has the additional feature that it leads to a natural interpretation of the resulting local time in terms of interval crossings. 

Let $V^1([0,T];\R)\subset D([0,T];\R)$ and $V^+([0,T];\R)\subset D([0,T];\R)$ be the space of all functions on $[0,T]$ with bounded variation (also called of finite total variation) and of all non-decreasing functions, respectively. Let us recall that for $[0,t] \subset [0,T]$ and $y\colon[0,T]\to\R$, the total variation of $y$ on the interval $[0,t]$ is given by
\begin{equation*}
  \TTV{y}{[0,t]}{} := \sup \bigg \{ \sum_{i=0}^{N-1} |y_{t_{i+1}}- y_{t_{i}}|\,:\, (t_i)_{i=0}^N\text{ is a partition of  } [0,t],\, N\in\N  \bigg\}.
\end{equation*}

\begin{definition}
  Given $x\in D([0,T];\R)$ and $\varepsilon>0$, a pair $(\phi^{\varepsilon}, -x^{\varepsilon})\in D([0,T];\R)\times V^1([0,T];\R)$ is called a \textup{solution to the Skorokhod problem} on $[-\varepsilon/2,\varepsilon/2]$ if the following conditions are satisfied:
  \begin{enumerate}
    \item[(i)] $x_t- x^{\varepsilon}_t =\phi^\varepsilon_t \in [-\varepsilon/2,\varepsilon/2]$ for every $t\in [0,T]$, 
    \item[(ii)] $x^{\varepsilon}= x^{\varepsilon\uparrow} - x^{\varepsilon\downarrow}$ with $x^{\varepsilon\uparrow}, x^{\varepsilon\downarrow}\in V^+([0,T];\R)\subset D([0,T];\R) $ and the corresponding measures  $\d x^{\varepsilon\uparrow}_t$  and $\d x^{\varepsilon\downarrow}_t $ are supported in $\{ t\in [0,T]\,:\, \phi^{\varepsilon}_t=\varepsilon/2\}$ and $\{ t\in [0,T]\,:\, \phi^{\varepsilon}_t=-\varepsilon/2\}$, respectively, 
    \item[(iii)] $\phi^{\varepsilon}_0= 0$.
  \end{enumerate}
\end{definition}

A solution to the above Skorokhod problem exists and is unique, see \cite[Proposition~2.7]{Lochowski2014}, and its properties are well studied in the literature, see, e.g., \cite{Kruk2007,Burdzy2009}. Let us emphasise that for any $\varepsilon>0$, $x^{\varepsilon}$ is a c{\`a}dl{\`a}g and piecewise monotonic path of bounded variation, which uniformly approximates $x$ with accuracy $\varepsilon/2$. 

While $f\circ y$ is of finite variation for all $y:[0,T]\to \R$ which are of finite variation if and only if $f$ is locally Lipschitz (see \cite{Leoni2017}), we can nonetheless assert that $f(x^{\varepsilon})$ is of finite variation for any $f\in W^{2,q}(\R)$, because $x^{\varepsilon}$ is a special function of finite variation: it is \emph{piecewise monotonic}, i.e., there exists a partition $0=a_0<a_1<\ldots <a_{N+1}=T$ of $[0, T]$ s.t. $x^{\varepsilon}$  is either increasing or decreasing  on each  $I_i$, where\footnote{\cite[Remarks~2.5 and~2.6]{Lochowski2014} imply that there is \emph{finite} number of such intervals: otherwise, the c{\`a}dl{\`a}g function $\phi^{\varepsilon}$ would have no left limit at the point $\lim_{i \rightarrow +\infty} a_i$, a contradiction.} 
$$
  I_i:=[a_i,a_{i+1}), i=0, \ldots, N-1, \quad I_N:=[a_N,a_{N+1}] ,
$$
see \cite[formula (2.4)]{Lochowski2014}, where even a more general Skorokhod problem is considered.
Thus,  keeping in mind integration by parts for the Lebesgue--Stieltjes integral, for $\varepsilon>0$ and $f\in W^{2,q}(\R)$, $q \ge 1$, we can define 
\begin{align}\label{def of intfxndx}
     \begin{split}
    \int_{0}^t f^\prime (x^{\varepsilon}_{s-}) \dd x_s
    & :=f^\prime (x^{\varepsilon}_{t})x_t-f^\prime (x^{\varepsilon}_{0})x_0-\int_{0}^t x_{s-}\dd f^\prime (x^{\varepsilon}_{s})
    -\sum_{0<s\le t}\Delta x_s\Delta f^\prime (x^{\varepsilon}_{s}),
    \end{split}
\end{align}
where $\int_{0}^t x_{s-}\dd f^\prime (x^{\varepsilon}_{s})$ exists as the Lebesgue--Stieltjes integral and we recall the convention $\int_{0}^t = \int_{(0,t]}$. For a brief summary of the theory of Lebesgue--Stieltjes integration, we refer to \cite[Chapter~4, Section~3.18]{RogersWilliams2000vol2}; we also remind the reader that, if $x,y$ are of finite variation, 
$$
  \int_0^t y_s \dd x_s= \int_0^t y_{s-} \dd x_s + \sum_{s\leq t} \Delta y_s \Delta x_s.
$$ 
We will define the pathwise local time as normalised limits of the numbers of interval crossings. To this end, for $x\in D([0,T];\mathbb{R})$, $z\in \R$, $\varepsilon>0$ and $t\in (0,T]$ we define the number of upcrossings by the path $x$ of the interval $(z-\varepsilon/2,z+\varepsilon/2)$ over the time $[0, t]$ by 
\begin{equation*}
  \Ucross{z,\varepsilon}x{[0,t]} := \sup_{n\in \N} \sup_{0\leq t_1<s_1<\cdots <t_n< s_n \leq t} \sum_{i=1}^{n} \1_{\{ x_{t_i}\leq z-\varepsilon/2\text{ and } x_{s_i} \ge z+\varepsilon/2\}}.
\end{equation*}
The number of downcrossings $\Dcross{z,\varepsilon}x{[0,t]}$ is defined analogously. We set
\begin{align}
\label{eq: def crossings}
  \cross{z,\varepsilon}x{[0,t]}:= \Dcross{z,\varepsilon}x{[0,t]} + \Ucross{z,\varepsilon}x{[0,t]}
\end{align} 
for the total number of crossings. 

\begin{definition}\label{def:local time via Skorokhod map}
  Consider a sequence $\rbr{c_n}_n$ such that $c_n>0$ and  $c_n \to 0$. For $x\in D([0,T];\mathbb{R})$ denote by $\rbr{\phi^n, -x^n}$ the solution to the Skorokhod problem on $[-c_n/2,c_n/2]$, $n\in \N$. 
  We denote $\mathbb{L}^{S}_p((c_n)_n)$ the set of all paths $x\in D([0,T];\mathbb{R})$ such that, for all $0< t\leq T$, 
    \begin{enumerate}
    \item[(i)] the sequence of functions 
      \begin{equation*}
        \R \ni z \mapsto  c_n \cdot \Cross{z,c_n}x{[0,t]},\quad n\in \N,
      \end{equation*}
      converges weakly in $L^p(\R)$ as $n \to \infty$; and
      
      \item[(ii)] the sequence of functions 
      \begin{equation*}
        \R \ni z \mapsto J_t(x^n, z) , \quad n\in \N,
      \end{equation*}
      defined by formula \eqref{eq: def J_t(x,u)}, converges weakly in $L^p(\R)$ to $J_t(x, \cdot)$ as $n \to \infty$.
  \end{enumerate}
A function $L \colon [0,T] \times \mathbb{R} \to \mathbb{R}$ which is the weak limit in $(ii)$ is called an $L^p$-\textup{interval crossing local time} of $x$ along $(c_n)_n$.
  \end{definition}

The corresponding pathwise Tananka--Meyer formula reads as follows.

\begin{proposition}\label{prop:Tananka-Meyer formula via Skorokhod}
  Suppose that $x\in\mathbb{L}^{S}_p((c_n)_n)$ for $p,q\geq 1$ with $1/p+1/q=1$. If $f\in W^{2,q}(\R)$, then the following limit exists and is finite
  \begin{equation*}
    \int_{0}^t f^\prime(x_{s-}) \dd x_s :=\lim_{n\to \infty}\int_{0}^t f^\prime (x^{n}_{s-}) \dd x_s, \quad t\in [0,T],
  \end{equation*}
  where the right-hand side is defined using~\eqref{def of intfxndx}, and the pathwise Tanaka--Meyer formula
   \begin{equation*}
    f(x_t)-f(x_0)  
    =\int_{0}^t f^\prime(x_{s-})\dd x_s+ \frac{1}{2} \int_{\R} L_t(x,u)  f^{\prime\prime}(\mathrm{d}u) +J_t^f(x),\quad t\in [0,T], 
  \end{equation*} holds with such definition of $\int_0^t f^\prime (x_{s-}) \dd x_s$ and with $J_t^f(x)$ as given in~\eqref{eq:sum of jumps}. 
\end{proposition}

Before proving Proposition~\ref{prop:Tananka-Meyer formula via Skorokhod}, we prove the following very intuitive lemma, where we write $g(x)$ for $g \circ x$. 

\begin{lemma}\label{property N}
  Let $I \subseteq \R$ be an open interval (i.e., $I$ is open and convex), $x\colon I\to [c,d]$ be c{\`a}dl{\`a}g and monotonic (i.e., increasing or decreasing), and $g\colon [c,d] \to \R$ be absolutely continuous and increasing. If $\d x$ is concentrated on a Borel set~$F$, then so is $\d g(x)$. 
\end{lemma}

\begin{proof}
  We can w.l.o.g. assume that $I=\R$, since otherwise we can trivially extend $x$ to $\R$ in a way that $\R\setminus I$ has $\d x$ mass $0$. We have to prove that the $\d x$ null set $E:=\R\setminus F$ is also a $\d g(x)$ null set. Denote with $\mathcal{L}$ the Lebesgue measure on $\R$. Let $y$ be c{\`a}dl{\`a}g monotonic, so if $I,J\subseteq \R$ are intervals with disjoint interiors, then so are $y(I),y(J)$ (even if $I \cap J =\emptyset$ does not imply $y(I) \cap y(J)= \emptyset$). Set $s(y):=1$ (resp. $=-1$) if $y$ is increasing (resp. decreasing). Since 
  \begin{align}\label{eq: image set}
    \int_\R \1_A(u) \dd y_u= s(y) \int_\R \1_{y(A)}(u) \dd u= s(y) \mathcal{L}(y(A))
  \end{align}
  holds (by definition of $\d y$) when $A$ is an interval, it holds whenever $A \subseteq \R$ is a countable union of intervals $(I_n)_n$ with disjoint interiors (because the interiors of $(y(I_n))_n$ are disjoint).

  Fix arbitrary $\epsilon>0$ and recall that there exists a $\delta>0$ s.t. $ \mathcal{L}(V)\leq \delta$ implies $\int \1_{V} \dd g=\mathcal{L}(g(V))\leq \epsilon$ whenever $V$ is a \emph{finite} union of intervals with disjoint interiors (by definition of absolute continuity), and thus whenever $V$ is a \emph{countable} union of intervals with disjoint interiors.

  Now cover $E$ with an open set $A \supseteq E$ s.t. $|\int \1_A \dd x|\leq \delta$; since $A$ is open, it is a countable union of disjoint open intervals, so \eqref{eq: image set} with $y=x$ gives $\mathcal{L}(x(A))\leq \delta$. Since $V:=x(A)\supseteq x(E)$ is a countable union of intervals with disjoint interiors we get $\mathcal{L}(g(x(A)))\leq \epsilon$, and so \eqref{eq: image set} with $y=g(x)$ gives $|\int \1_{A} \dd g(x)|\leq \epsilon$. Thus $| \int \1_{E} \dd g(x)|\leq \epsilon$ for any $\epsilon>0$, concluding the proof.
\end{proof} 

\begin{proof}[Proof of Proposition~\ref{prop:Tananka-Meyer formula via Skorokhod}]
  We introduce first slightly modified numbers of interval (up-, down-) crossings by replacing $\leq, \geq $ with $<, >$ in the inequality involving $x_{t_i}$ in the definition of up- and down- crossings: for $z\in \R$, $\varepsilon\ge0$, $t\in (0,T]$ and $x\in D([0,T];\R)$ we set
  \begin{align*}
    &\Ucrosstilde{z,\varepsilon}x{[0,t]} := \sup_{n\in \N} \sup_{0\leq t_1<s_1<\cdots <t_n< s_n \leq t} \sum_{i=1}^{n} \1_{\{ x_{t_i}< z-\varepsilon/2\text{ and } x_{s_i} \ge z+\varepsilon/2\}},\\
    &\Dcrosstilde{z,\varepsilon}x{[0,t]} := \sup_{n\in \N} \sup_{0\leq t_1<s_1<\cdots <t_n< s_n \leq t} \sum_{i=1}^{n} \1_{\{ x_{t_i}> z+\varepsilon/2\text{ and } x_{s_i} \le z+\varepsilon/2\}},\\
    & \crosstilde{z,\varepsilon}x{[0,t]}:= \Dcrosstilde{z,\varepsilon}x{[0,t]} + \Ucrosstilde{z,\varepsilon}x{[0,t]}.
  \end{align*}
  As $f\in W^{2,q}(\R)$, $f^\prime$ can be decomposed as the difference of two increasing AC (Absolutely Continuous) functions; since the result we want to prove is linear in $f$, we can assume w.l.o.g. that $f^\prime$ is increasing and AC. Moreover, since $x$ (as defined on $[0,t]$) is  bounded, and the result only depends on the behaviour of $f$ on $[\inf x, \sup x]$, we can additionally assume w.l.o.g. that $f^{\prime\prime}$ has compact support. As the proposition holds trivially for affine functions, thanks to \eqref{eq:representation of f} we may further assume that
  \begin{equation*}
    f(x) =  (|\cdot | * f^{\prime\prime})(x), \quad x\in \R .
  \end{equation*}
  
  Let us consider the integral $\int_{0}^t f^\prime (x^{n}_{s-}) \dd x_s$. For $t\in [0,T]$ we have 
  \begin{align}\label{eq:jeden}
    \begin{split}
    \int_{0}^t f^\prime (x^{n}_{s-}) \dd x_s
    & =f^\prime (x^{n}_{t})x_t-f^\prime (x^{n}_{0})x_0-\int_{0}^t x_{s-}\dd f^\prime (x^{n}_{s})
    -\sum_{0<s\le t}\Delta x_s\Delta f^\prime (x^{n}_{s}),
    \end{split}
  \end{align}
  where $\int_{0}^t x_{s-}\dd f^\prime (x^{n}_{s})$ is the Lebesgue--Stieltjes integral. Further, we have 
  \begin{align} \label{eq:dwa}
    \begin{split}
     \int_{0}^t x_{s-}\dd f^\prime (x^{n}_{s})+\sum_{0<s\le t}\Delta x_s \Delta f^\prime (x^{n}_{s})
    & =\int_{0}^t x_{s}\dd f^\prime (x^{n}_{s}) \\
    & =\int_{0}^t (x_{s} - x_s^n) \dd f^\prime (x^{n}_{s})+\int_{0}^t x_s^n \dd f^\prime (x^{n}_{s}).
    \end{split}
  \end{align}
  To calculate the first integral we use the properties of $f'$ and $x^n$.   

  Recall that the positive (resp. negative) part of $\dd x^n$ is concentrated on $\{x-x^n= c_n/2\}$ (resp. $\{x-x^n=-c_n/2\}$). Thus, the identity 
  \begin{equation}\label{eq:ttv}
    \int_{I} (x_{s} - x_s^n) \dd f^\prime (x^{n}_{s})=\frac{c_n}{2}\TTV{f^\prime \rbr{x^{n}_{\cdot}}}{I}{},
  \end{equation} 
  holds if $I$ is the interior of an interval on which $x^n$ is increasing (resp. decreasing), by Lemma~\ref{property N}, and if $I$ is a singleton, since in that case it reduces to the identity 
  $$
    (x_{s} - x_s^n)  \Delta f^\prime (x^{n}_{s})=\frac{c_n}{2}|\Delta f^\prime (x^{n}_{s})|.
  $$
  Since $x^n$ is piecewise monotonic, we conclude that \eqref{eq:ttv} holds for $I=[0,t]$.
  
  Using~\eqref{eq:jeden}, \eqref{eq:dwa} and \eqref{eq:ttv}, we finally arrive at 
  \begin{align}\label{eq:formula_gen}
    \int_{0}^t f^\prime (x^{n}_{s-}) \dd x_s
     =f^\prime (x^{n}_{t})x_t-f^\prime (x^{n}_{0})x_0-\int_{0}^t x^n_{s}\dd f^\prime (x^{n}_{s})
     -\frac{c_n}{2}\TTV{f^\prime \rbr{x^{n}_{\cdot}} }{[0,t]}{}.
  \end{align}
  Let us note that the right side of~\eqref{eq:ttv} may be also calculated using  the following generalisation of the Banach indicatrix theorem:
  \begin{equation} \label{eq:Ban_ind}
    \TTV{f^\prime \rbr{x^{n}_{\cdot}}}{[0,t]}{}=\int_{\R}N^{y}\left(f^\prime \rbr{x^{n}_{\cdot}},[0,t]\right)\dd y ,
  \end{equation}
  where $N^{y}\left(g,[0,t]\right)$ is the number of up- and down- crossings of the level $y$ by c{\`a}dl{\`a}g $g$, as defined in \cite[Remark~1.3]{Lochowski2017}, which is closely related to the number of crossings $n^{y,\varepsilon}$  via the relation
  \begin{equation*}
    N^{y}\left(g,[0,t]\right) = \lim_{\varepsilon \to 0+}   \cross{y,\varepsilon}g{[0,t]},
  \end{equation*}
  of which we will not make use, and which can be proved similarly to \cite[Remark~1.4]{Lochowski2017}. Moreover, the relationship
  \begin{align}\label{eq:Ban_ind2}
    \int_{\R}N^{y}\left(f^\prime \rbr{h_{\cdot}},[0,t]\right)\dd y = \int_{\R}N^{z}\left(h,[0,t]\right)\dd{f^\prime\left(z\right)}
  \end{align}
  clearly holds for any monotonic $h\colon  I\to \R$ defined on an open interval $I$, and thus holds for any completely monotonic $h$.
  
  Thus, equation~\eqref{eq:formula_gen} may be rewritten as 
  \begin{align*}
    \int_{0}^t x^n_{s}\dd f^\prime (x^{n}_{s})
     =f^\prime (x^{n}_{t})x_t-f^\prime (x^{n}_{0})x_0-\int_{0}^t f^\prime (x^{n}_{s-}) \dd x_s
     -\frac{c_n}{2}\TTV{f^\prime \rbr{x^{n}_{\cdot}} }{[0,t]}{}
  \end{align*}
  which, thanks to \eqref{eq:Ban_ind}, \eqref{eq:Ban_ind2}, takes the form
  \begin{align} \label{eq:cztery}
    \int_{0}^t x^n_{s}\dd f^\prime (x^{n}_{s})
     =f^\prime (x^{n}_{t})x_t-f^\prime (x^{n}_{0})x_0-\int_{0}^t f^\prime (x^{n}_{s-}) \dd x_s
     -\frac{c_n}{2}\int_{\R}N^{z}\left(x^{n},[0,t]\right)\dd{f^\prime\left(z\right)}.
  \end{align}
  Now we will compute an alternative expression for  $\int_{0}^t x^n_{s}\dd f^\prime (x^{n}_{s})$. Since $x^n$ and $f^\prime \rbr{x^{n}_{\cdot}}$ have finite total variation the rules of the Lebesgue--Stieltjes integral (integration by parts and the substitution rule) apply here and we have
  \begin{equation}\label{eq:trzyy}
    \int_{0}^t x^n_{s}\dd f^\prime (x^{n}_{s}) 
    =f^\prime (x^{n}_{t})x^n_t-f^\prime (x^{n}_{0})x^n_0-\int_{0}^t f^\prime (x^{n}_{s}) \dd x^n_{s} 
    +\sum_{0<s\le t}\Delta x^n_s\Delta f^\prime (x^{n}_{s}) 
  \end{equation}
  and 
  \begin{equation}\label{eq:trzy}
    \int_{0}^t f^\prime (x^{n}_{s}) \dd x^n_{s}  = f(x^{n}_{t}) - f (x^{n}_{0}) -\sum_{0<s\le t}\big ( \Delta f(x^n_{s})-f^\prime (x^n_{s} )\Delta x^n_s \big).
  \end{equation}
  Since
  $$
    \Big(\Delta f(x^n_{s})-f^\prime (x^n_{s} )\Delta x^n_s \Big) + \Delta x^n_s\Delta f^\prime (x^{n}_{s})= \Delta f(x^n_{s})-f^\prime (x^n_{s-} )\Delta x^n_s ,
  $$
  whose sum over $s\leq t$ equals $J^f_t(x^n)$, substituting in \eqref{eq:trzyy} the value for $\int_{0}^t f^\prime (x^{n}_{s}) \dd x^n_{s}$ obtained from \eqref{eq:trzy} we get 
  \begin{align}\label{eq:piec}
    \int_{0}^t x^n_{s}\dd f^\prime (x^{n}_{s})   &  =  f^\prime (x^{n}_{t})x^n_t-f^\prime (x^{n}_{0})x^n_0-  ( f(x^{n}_{t}) - f (x^{n}_{0})) +J_t^f(x^n).
  \end{align}
  Finally, equating the RHS of \eqref{eq:piec} and \eqref{eq:cztery}  we get 
  \begin{align}\label{eq:formula_gen_1}
    \begin{split}
    f(x^{n}_{t}) - f (x^{n}_{0}) 
    & = \int_{0}^t f^\prime (x^{n}_{s-}) \dd x_s + \frac{c_n}{2} \int_{\R}N^{z}\left(x^{n},[0,t]\right)\dd{f^\prime\left(z\right)} +J_t^f(x^n) \\
    & \quad  - f^\prime (x^{n}_{t})\rbr{x_t - x^n_t} + f^\prime (x^{n}_{0})\rbr{x_0 - x^n_0}.
    \end{split}
  \end{align}
  Let us now compute 
  $$
    \lim_{n\to \infty} c_n \int_{\R}N^{z}\left(x^{n},[0,t]\right)\dd f^\prime\left(z\right) .
  $$
  Since the set of local extrema (maxima and minima) of any function $f\colon\R \to \R$ is countable (see e.g. \cite[Lemma~5.1]{PoseyVaughan1986}), the numbers $N^{z}\left(x^n,[0,t]\right)$ and $\crosstilde{z,0}{x^n}{[0,t]}$ are equal for all $z\in\R$ except a countable set, because they are equal if $z\notin \cbr{x_0, x_t}$ and $z$ is not a local extremum of $x$. Similarly, for all $z\in\R$ except a countable set, the numbers $\crosstilde{z,c_n}{x}{[0,t]}$ and $\cross{z,c_n}{x}{[0,t]}$ are equal, because they are equal if $z\notin \cbr{x(0)\pm c_n/2, x(t)\pm c_n/2}$ and $z\pm c_n/2$ are not local extrema of $x$. Next, by \cite[Lemma~3.3 and~3.4]{Lochowski2014}, $\crosstilde{z,0}{x^n}{[0,t]}$ and $ \crosstilde{z,c_n}x{[0,t]}$ differ by at most $2$. Thus $N^{z}\left(x^n,[0,t]\right)$ and $\cross{z,c_n}{x}{[0,t]}$ differ by at most $2$ for all but a countable number of  $z\in\R$. Using this observation and noticing that  $N^{z}\left(x^n,[0,t]\right)=\cross{z,c_n}x{[0,t]}=0$ when $z<\inf_{s\in[0,t]}x_s-c_n/2$ or $z>\sup_{s\in[0,t]}x_s+c_n/2$ we have that 
  \begin{align*}
    \lim_{n\to \infty} c_n \int_{\R}N^{z}\left(x^{n},[0,t]\right)\dd f^\prime\left(z\right)  
    =  \lim_{n\to \infty} \int_{\R}c_n\cdot\cross{z,c_n}x{[0,t]}f^{\prime\prime} \left(z\right)\dd{z}  
    = \int_{\R}L_t(z) f^{\prime\prime} \left(z\right)\dd z,  
  \end{align*}
  where the last equality follows from the first assumption in Definition~\ref{def:local time via Skorokhod map}. Also, by the second assumption in Definition~\ref{def:local time via Skorokhod map} and Lemma~\ref{JfJut} 
  \begin{align*}
    \lim_{n\to \infty} J_t^f(x^n)  = \lim_{n\to \infty} \int_{\R} J_t (x^n, y) f^{\prime\prime}(y) \dd y  
    = \int_{\R} J_t (x, y) f^{\prime\prime}(y) \dd y = J_t^f(x). 
  \end{align*}
  The last two limits together with~\eqref{eq:formula_gen_1} give the thesis.
\end{proof}

\begin{remark}\label{condition_jumps}
  To apply Proposition~\ref{prop:Tananka-Meyer formula via Skorokhod} we need to know when $J_t (x^n, \cdot)$ converge weakly in $L^p\rbr{\R}$ to $J_t (x, \cdot) \in L^p\rbr{\R}$ and 
  $c_n\cdot\cross{\cdot,c_n}x{[0,t]}$ converges weakly in $L^p\rbr{\R}$ to some $L_t \in L^p\rbr{\R}$. However, in general, it is not even clear when $c_n\cdot\cross{\cdot,c_n}x{[0,t]}$ and $J_t (x^n, \cdot)$ belong to $L^p\rbr{\R}$ although we now give some sufficient criteria.  
  If for some $r>0$, the $r$-variation is finite, i.e., 
  \begin{equation*}
    V^r\rbr{x,[0,T]} := \sup \bigg \{ \sum_{i=0}^{N-1} |x_{t_{i+1}}- x_{t_{i}}|^r\,:\, (t_i)_{i=0}^N\text{ is a partition of  } [0,T],\, N\in\N  \bigg\}<\infty,
  \end{equation*}
  then $c_n \cdot \cross{\cdot,c_n}x{[0,t]}$ is bounded (and is equal $0$ outside a compact subset of $\R$) and thus belongs to  $L^p\rbr{\R}$ for all $t\leq T$. It follows from the easy estimate: for any $z \in \R$
  \[
    \cross{z,c_n}x{[0,t]} \le \frac{V^r\rbr{x,[0,t]}}{c_n^r}.
  \] 
  Unfortunately, this observation does not yield any condition which guarantees $L_t \in L^p\rbr{\R}$, except in the rather trivial case $r \leq 1$.
  
  Similarly as in Remark~\ref{Jfinite} we have that if  $p \in [1,\infty)$ and $\sum_{0<s\leq t} |\Delta x_s|^{1+1/p} < \infty$ then $J_t (x^n, \cdot), J_t (x, \cdot) \in L^p\rbr{\R}$: this follows from Minkowski's inequality and the fact that for any $s>0$, $\left|\Delta x^n_{s} \right| \le \left|\Delta x_{s} \right|$ (see \cite[(2.5)]{Lochowski2014} or \cite[Section~2]{Lochowski2014b}). 

  Since $x^n \rightarrow x$ \emph{uniformly}, there exists $c\in \R$ s.t., for all $s\in [0,t]$,  
  \begin{align}\label{eq: xn to x}
  \begin{split}
    &|x^n_s|\leq c \text{ for all } n\in \N ,  \\
    &\exists \lim_n |x^n_s-u| \1_{  \llbracket x^n_{s-}, x^n_{s} \rrparenthesis }(u)
    = |x_s-u| \1_{  \llbracket x_{s-}, x_{s} \rrparenthesis }(u) \text{ for all } u \neq x_{s-}, x_{s} .
  \end{split}  
  \end{align}
  Now let $q$ be s.t. $1/p+1/q =1$, and fix any $f \in W^{2,q}\rbr{\R}$ and $s\in [0,t]$. Since $f''$ is locally integrable, it follows from~\eqref{eq: xn to x} and the dominated convergence theorem that 
  $$
    \int_{\R} |x^n_s-u| \1_{  \llbracket x^n_{s-}, x^n_{s} \rrparenthesis }(u)f^{\prime \prime} (u) \dd u \to \int_{\R} |x_s-u| \1_{  \llbracket x_{s-}, x_{s} \rrparenthesis }(u)f^{\prime \prime} (u) \dd u \quad \text{ as } n \to \infty .
  $$
  We can then apply again the dominated convergence theorem to obtain weak convergence of $J_t (x^n, \cdot)$ to $J_t (x, \cdot)$ in $L^p\rbr{\R}$, using the domination 
  \begin{equation*}
    \int_{\R} |x^n_s-u| \1_{  \llbracket x^n_{s-}, x^n_{s} \rrparenthesis }(u)|f^{\prime \prime} (u)| \dd u  \le C_p \left|\Delta x^n_{s} \right|^{1+\frac{1}{p}} \Vert f^{\prime \prime}\Vert_{L^q} \le C_p \left|\Delta x_{s} \right|^{1+\frac{1}{p}} \Vert f^{\prime \prime}\Vert_{L^q},
  \end{equation*}
  which follows from the estimate $\left|\Delta x^n_{s} \right| \le \left|\Delta x_{s} \right|$, H{\"o}lder's inequality and \eqref{eq:int |b-u|^p}.
\end{remark}

\section{Construction of local times for c{\`a}dl{\`a}g semimartingales}\label{sec:construction of local limes}

The purpose of this section is to give probabilistic constructions of the pathwise local time, as introduced in Definitions~\ref{def:L^pLT}, \ref{def:Wuermli local time} and \ref{def:local time via Skorokhod map}, for c{\`a}dl{\`a}g semimartingales. In particular, we show that all three definitions agree a.s.\ and coincide with the classical probabilistic notion of local times for c{\`a}dl{\`a}g semimartingales.

\subsection{Local times via discretisation and as occupation measure}

Let $(\Omega, \mathcal{F},\mathbb{F},\P)$ be a filtered probability space where the filtration $\mathbb{F}:=(\mathbb{F}_t)_{t\in [0,\infty)}$ is supposed to satisfy the usual conditions. Given a c{\`a}dl{\`a}g semimartingale $X=(X_t)_{t\in [0,\infty)}$ and $u \in \R$, one can define $ \mathcal{J}_t(u)(\omega):= J_t(X_\cdot(\omega),u)$, with $J_t(x,u)$ given by~\eqref{eq: def J_t(x,u)}, and the increasing c{\`a}dl{\`a}g adapted process $\bL(u)$ by
\begin{equation}\label{eq:stoc tanaka}
  2\bL_t(u) := |X_t-u|-  |X_0-u| -\int_{(0,t]} \textup{sign}(X_{s-}- u) \dd X_s .
\end{equation}
It can then be shown that there exists a jointly measurable version of $ \bL_t(u,\omega)$ such that the family of processes $\mathcal{L}=2\bL-2\mathcal{J}$, called the (classical) local time of $X$, satisfies the Tanaka--Meyer formula~\eqref{eq:Tanaka formula} for $x=X(\omega)$ $\P(\d \omega)$-a.e., is c{\`a}dl{\`a}g in $t$ and is jointly measurable: see\footnote{Recall the identity~\eqref{eq:Two expressions for J} and notice that the notations used in \cite{Protter2004} differ from ours: he calls $A^a$ what we call $2\bL(u)$.} \cite[Chapter~4, Section~7]{Protter2004}. 

In the following we denote by $L^p(\mu)$ the $L^p$-space with respect to a measure $\mu$. If $\pi=(\tau_k)_{k\in \N}$, where $\tau_k$ are $[0,\infty]$-valued random variables such that $ \tau_0=0, \tau_k\leq\tau_{k+1}$ with $\tau_k<\tau_{k+1}$ on $\{\tau_{k+1}<\infty \}$, and $\lim_{k\to \infty}\tau_{k}=\infty$, then $\pi$ is called a random partition. If moreover $\{\tau_k\leq t \} \in \mathcal{F}_t$ for all $k,t$, then $\pi$ is called an optional partition. We recall $K_s^{\pi}$ was defined in~\eqref{eq:Kpi}. The following is the main theorem of this subsection.  

\begin{theorem}\label{GenExistsLT}
  Assume that $f\colon \R\to \R$ is a difference of two convex functions, $(\pi^n)_n$ are optional partitions of $[0,\infty)$ such that $|\pi^n \cap [0,t] |\to 0$ a.s. for all $t$ and $X=(X_t)_{t\in [0,\infty)}$ is a c{\`a}dl{\`a}g semimartingale. Then, there exists a subsequence $(n_k)_k$ such that, for $\omega$ outside of a $\P$-null set (which may depend on $f''$), 
  \begin{align*}
    \sup_{s\in [0,t]} \left| K_s^{\pi^{n_k}(\omega)}(X_\cdot(\omega),u) - \bL_{s}(\omega,u) \right|  \to 0 \quad  \text{ in } L^p(|f''|(\d u)) \text{ as } k\to \infty    
  \end{align*} 
  simultaneously for all $p\in[1,\infty)$, $t\in [0,\infty)$.
\end{theorem} 

\begin{remark}
  Theorem~\ref{GenExistsLT} says that the pathwise crossing time $K_\cdot^{\pi^{n}}(X_\cdot,u)$ sampled along optional partitions $(\pi^n)_n$ (defined applying~\eqref{eq:Kpi} to each path $X_\cdot(\omega)$ and partition $\pi^n(\omega)$) converges to $\bL(u)$. Applying Theorem~\ref{GenExistsLT} with $f(x)=x^2/2$ gives in particular that $\P(\d\omega)$-a.e.\ $X(\omega)\in \mathbb{L}_p^W((\pi^{n_k})_k)\subset \mathbb{L}_p((\pi^{n_k})_k)$ for all $p<\infty$ and $T>0$, i.e., the $L^p$-level crossing time and the occupation local time exist for a.e. paths of a semimartingale. Indeed, $K_t^{\pi^{n_k}}(X,\cdot) \to \bL_{t}(\cdot) $ strongly (and thus weakly) in $L^p(\R)$ for a.e. $\omega$, locally uniformly in $t$. 
\end{remark}

To prove the previous theorem we need some preliminaries. Given $p\in [1,\infty)$ we denote by $\cS^p$ the set of c{\`a}dl{\`a}g special semimartingales~$X$ which satisfy
\begin{equation*}
  \|X\|_{\cS^p}:=\left\|[Y]_{\infty}^{1/2} \right\|_{L^p(\P)} + \left\| \int_0^\infty\ \dd|V|_t \right\|_{L^p(\P)}<\infty ,
\end{equation*}
where $X=Y+V$ is the canonical semimartingale decomposition of $X$, 
$$
  [Y]_t:=Y_t^2-2 \int_0^{t} Y_{s-} \dd Y_s 
$$
is the quadratic variation of the martingale $Y$, and $|V|_t$ is the variation up to time $t$ of the predictable process $(V_t)_{t\in [0,\infty)}$. We recall the existence of $c_p<\infty$ such that the inequality 
\begin{align}\label{GenBDG}
  \bigg\|\sup_{t\in [0,\infty)} \left| X_t \right|\bigg \|_{L^p(\P)}  \leq c_p \|X\|_{\cS^p} \, ,
\end{align} 
holds for all local martingales $X$ (this being one side of the Burkholder--Davis--Gundy inequalities), and thus also trivially extends to all $X\in \cS^p$. The core of Theorem~\ref{GenExistsLT} is the following more technical statement.

\begin{proposition}\label{ExistsLT}
  Let $(\pi^n)_n$ be optional partitions of $[0,\infty)$ such that $|\pi^n \cap [0,t] |\to 0$ a.s. for all $t$. If $X\in \cS^p$ for $p\in[1,\infty)$, and 
  \begin{align*}
    h^{\pi^n}(u):= \left\| \sup_{t\in [0,\infty)} \left| K_t^{\pi^n(\omega)}(X_\cdot(\omega),u) - \bL_{t}(\omega,u) \right| \right\|_{L^p(\P)}, \quad u\in \R ,
  \end{align*} 
  then, for every $u\in \R$, $h^{\pi^n}(u)\to 0$ as $n\to \infty$ and 
  $0\leq  h^{\pi^n}(u)\leq c_p \|  X \|_{\cS^p} $ for all $n\in \N$.
\end{proposition} 

As discussed in detail in \cite{Davis2018} after Theorem~6.2, for a continuous process $X$ and properly chosen $(\pi^{n})_n$ the convergence of $K_\cdot^{\pi^{n}}(X_\cdot,u)$ is closely related to the number of upcrossings of $X$ from the level $u$ to the level $u+\varepsilon_n > u$. While stronger versions of the above theorems have already appeared in the case of continuous semimartingales (the strongest being \cite[Theorem~II.2.4]{Lemieux1983}), in the c{\`a}dl{\`a}g setting we were only able to locate in the literature a version of Theorem~\ref{GenExistsLT} where, under the strong assumption that $\sum_{s\leq t} |\Delta X_s|<\infty$ a.s., the $L^p(|f''|(\d u))$ convergence is replaced by pointwise convergence for all but countably many values of $u$, see \cite[Theorem~III.3.3]{Lemieux1983}. Thus, compared to the literature, our method provides a novel strong conclusion, with the benefit of a simple proof. Other differences are that we consider the crossing time instead of the number of upcrossings, and we use any optional partitions such that $|\pi^n|\to 0$ instead of ``Lebesgue partitions'' (in the language of \cite{Davis2018}).

\begin{proof}[Proof of Proposition~\ref{ExistsLT}]
  Consider the convex function $f(x):=|x-u |$ and let us take its left-derivative $\sign(x-u)$ and its second (distributional) derivative $2\delta_{u}$. Subtracting from the discrete-time Tanaka--Meyer formula~\eqref{eq:discrete tanaka-meyer formula} its continuous-time stochastic counterpart~\eqref{eq:stoc tanaka} and considering the process $K_t^{\pi^n}(u)(\omega):=K_t^{\pi^n}(X_\cdot(\omega),u)$ we obtain 
  \begin{align}\label{LocTDifInt}
    0= \int_0^t  (H_s^{\pi^n}(u) -H_s(u)) \dd X_s+ 2( K_t^{\pi^n}(u)-\bL_t(u)) ,
  \end{align} 
  where for $\pi^n=(\tau_i^n)_i$ by $H^{\pi^n}$ and $H(u)$ we denote the predictable processes 
  \[
    H_s^{\pi^n}(u):=\sum_{i}\sign(X_{\tau_i^n}-u) \1_{(\tau_i^n, \tau_{i+1}^n] }(s) \quad \text{and}\quad  H_s(u):=\sign(X_{s-}-u) .
  \]
  Now $h^{\pi^n}(u)\to 0$ for each $u\in \R$ follows from~\eqref{GenBDG} and~\eqref{LocTDifInt} if we show that $\int_0^{\cdot} H_s^{\pi^n}(u) \dd X_s \to \int_0^{\cdot} H_s(u) \dd X_s$ in $\cS^p$. To this end fix $n$ and $u$ and notice that from 
  $$
    H_s^{\pi^n}(u)=\sign (X_{\tau_i^n}-u),\quad \text{for $i$ such that } \tau_i^n<s  \leq  \tau_{i+1}^n 
  $$
  and $|\pi^n \cap [0,t] |\to 0$ a.s. for all $t$ it follows that $H_s^{\pi^n}(u)\to H_s(u)$ a.s. for all $s$. Since $|H_s^{\pi^n}(u) -H_s(u)|\leq 2$, it follows that $\int_0^{\cdot} H_s^{\pi^n}(u)  \dd X_s \to \int_0^{\cdot} H_s(u)  \dd X_s$ in $\cS^p$  (by the dominated convergence theorem) and that 
  \[
    h^{\pi^n}(u) \leq  \frac{c_p}{2} \left\|  \int_0^{\cdot}  (H_s^{\pi^n}(u) -H_s(u))\dd X_s  \right\|_{\cS^p}  \leq c_p \|  X \|_{\cS^p}  \quad  \text{ for all   }  u\in \R  , \,
  \]
  concluding the proof.
\end{proof} 

\begin{proof}[Proof of Theorem~\ref{GenExistsLT}]
  Let $(\tau_m)_m$ be a sequence of stopping times which prelocalizes $X$ to $\cS^p$ (see e.g. \cite[Chapter~V, Theorem~4]{Protter2004}), i.e., $\tau_m \uparrow \infty$ a.s. and $X^{\tau_m -}\in \cS^p$ for all $m$. Let $\mu_i(A):=|f''|(A\cap [-i,i])$  and set, for $T>0$,
  \begin{equation*}
    G_n:= \textstyle G_n(\omega,T,u):=   \sup_{t\leq T} | K_t^{\pi^{n_k}(\omega)}(X_\cdot (\omega),u) - \bL_{t}(u,\omega) |
  \end{equation*}
  and $G_n^m:= \1_{\{T< \tau_m \}} G_n$. Since $\mu_i$ is a finite measure, Proposition~\ref{ExistsLT} implies that, as $n\to \infty$, $G_n^m$ converges  to  $0$ in $L^p(\P\times \mu_i)$  for all $m,i\in \N$ and $T\geq 0$. By Fubini's theorem  $||G_n^m||_{L^p(\mu_i)}$ converges to zero in $L^p(\P)$, and so passing to a subsequence (without relabelling) we find that, for every $\omega$ outside a $\P$-null set $N_{i,m}^{p,T}$,  $||G_n^m(\omega,T,\cdot)||_{L^p(\mu_i)}\to 0$. Then along a diagonal subsequence we obtain that $G_n^m(\omega,T,\cdot)$ converges to $0$ in $L^p( \mu_i)$ for all $i,m,p,T \in \N \setminus \{0\}$ for every $\omega$ outside the null set $N_{f''}:=\cup_{i,m,T,p\in \N\setminus \{0\}} N_{i,m}^{p,T}$. Since $G_n=G_n^m$ on $\{T< \tau_m \}$, $G_n\to 0$ in $L^p( \mu_i)$ for all $i,p,T \in \N \setminus \{0\}$ for every $\omega$ outside  $N_{f''}$. Since outside a compact set $G_n(\omega,T,\cdot)=0$ for all $n$, convergence in $L^p(\mu_i)$ for arbitrarily big $i,p$ implies convergence in $L^p(|f''|(\d u))$ for all $p\in [1,\infty)$. Since $G_n(\omega,\cdot,u)=0$ is increasing, convergence for arbitrarily big $T$ implies convergence for all $T\in[0,\infty)$.
\end{proof} 

\subsection{Local times via interval crossings}

Recall the definition of $L^p$-interval crossing local time of a deterministic path along a sequence of positive reals tending to $0$ in Definition~\ref{def:local time via Skorokhod map}. In this subsection we prove the following theorem.   

\begin{theorem}\label{GenExistsLTlevelcross}
  Let $X=(X_t)_{t\in [0,\infty)}$ be a c{\`a}dl{\`a}g semimartingale and $T >0$. There exist a $\P$-null set $E$ such that for any $\omega \in \Omega \setminus E$ and any sequence of positive reals $\rbr{c_n}_{n\in \N}$ which converges to $0$, $x_t = X_t\rbr{\omega}$, $t\in [0,T]$, belongs to $\mathbb{L}^{S}_1((c_n)_{n \in \N})$ and for any $t \in [0, T]$ the $L^1$-interval crossing local time of $x$ along $\rbr{c_n}_{n \in \N}$, $L_t$, coincides (in $L^{1}\rbr{\R}$) with the classical local time of $X$, $\mathcal{L}_t$. 
\end{theorem}

We note a difference in the above result when compared with Theorem \ref{GenExistsLT}. In the former, we obtained pathwise convergence on a subsequence and outside a null set which depended on the discretisation, i.e., on the optional partitions $(\pi^n)_n$ of $[0,\infty)$. Here, the method of discretisation is fixed and implicit in the Skorokhod problem, however we are able to obtain pathwise convergence, outside of a common null set $E$, simultaneously for all sequences $(c_n)$.

As noted already after the statement of Proposition~\ref{ExistsLT}, a similar result was proven in \cite[Theorem~III.3.3]{Lemieux1983}, namely that for any c{\`a}dl{\`a}g semimartingale $X$, as $c\to 0$, $c \cdot \cross{u,c}{X}{[0,t]}\to \mathcal{L}_t^u$ a.s. for all but countably many $u\in \R$, where $\text{n}^{u,c}$  was defined in \eqref{eq: def crossings}. However this was only established for semimartingales whose jumps are a.s. summable, i.e., $\sum_{0<s\le t} \left| \Delta X_s\right| < \infty$ for any $t>0$.

Theorem~\ref{GenExistsLTlevelcross} is easily proved using the following technical statement (of independent interest), about the quantity
\begin{align*}
  Q_t^{z,d}:= d\cdot\cross{z,d}X{\left[0,t\right]}-\frac{1}{d}\int_{z-d/2}^{z+d/2}{\mathcal L}_{t}^{u}\dd u,\quad t\in [0,\infty).
\end{align*} 

\begin{theorem}\label{conv}
  Let $X=(X_t)_{t\in [0,\infty)}$ be a c{\`a}dl{\`a}g semimartingale and $\mathcal{L}_t^u$, $t\geq 0$, $u\in \R$, its local times. If $(d_k)_{k\in \N}$ is a sequence of positive reals such that $\sum_{k\in \N} d_k < \infty$ then 
  \begin{align}\label{eq:conv} 
    \int_\R |Q_t^{z,d_k}| \dd z \to 0  \quad \P(\d \omega)\text{-a.e. as }k\to \infty ,
  \end{align} 
  and if $X \in \cS^{2p}$ for $p\in [1,\infty)$ and $|X|$ is bounded by a constant then  for any $t\in [0,\infty)$
  \begin{align}\label{eq: conv S2}
    \exists \lim_{d \downarrow 0}  \left\Vert \int_\R |Q_t^{z,d}| \dd z \right\Vert_{L^p(\P)}= 0.
  \end{align} 
\end{theorem}

Let us now prove Theorem~\ref{GenExistsLTlevelcross}; the rest of the subsection will be devoted to the proof of Theorem~\ref{conv}.

\begin{proof}[Proof of Theorem~\ref{GenExistsLTlevelcross}]
  By standard properties of convolutions, for example \cite[Theorem~8.14]{Folland1999}, $\frac{1}{d}\int_{z-d/2}^{z+d/2}{\mathcal L}_{t}^{y}\dd y \to \mathcal L_t^z$ in $L^1\rbr{\R}$ as $d\to 0+$. By this and Theorem~\ref{conv} there exists a $\P$-null set $E_1$ such that for any $\omega \in \Omega_1 = \Omega \setminus E_1$ and $x = X\rbr{\omega}$ the limit of $d_k \cdot \cross{\cdot,d_k}{x}{[0,t]}$ in $L^1\rbr{\R}$ (thus also the weak limit in $L^1\rbr{\R}$), where for example $d_k = k^{-2}$, exists and is equal $\mathcal{L}_t(\cdot) $ as $k \to \infty$. Now, for the given sequence $\rbr{c_n}_n$ and  $n$ such that $c_n \le 1/2$ define $k(n)$ to be such natural number that $d_{k(n)+1} < c_n \le d_{k(n)}$. For such $n$ we have bounds
  \begin{align}\label{bound} 
  \begin{split}
    \left(\frac{d_{k(n) + 1} }{d_{k(n)} }\right) d_{k(n)} \cdot \cross{\cdot ,d_{k(n)}}x{\left[0,t\right]}
    &\le c_n \cdot \cross{\cdot ,c_n }x{\left[0,t\right]} \\
    &\le \left(\frac{d_{k(n)} }{d_{k(n)+1} }\right) d_{k(n)+1} \cdot \cross{\cdot ,d_{k(n)+1}}x{\left[0,t\right]}.
  \end{split}
  \end{align}
  Notice, that since $d_k/d_{k+1} \to 1$ as $k \to \infty$, we have that for any $\omega \in \Omega_1$ the limits in $L^1\rbr{\R}$ of both -- lower and upper -- bounds in \eqref{bound} as $n \to \infty$ coincide with the limit of $d_k \cdot \cross{\cdot,d_k}{x}{[0,t]}$ which is equal $\mathcal{L}_t(\cdot) $. Thus for $\omega \in \Omega_1$,  $c_n \cdot \cross{\cdot,c_n}{x}{[0,t]}$ tends in $L^1\rbr{\R}$ to the same limit $\mathcal{L}_t(\cdot) $.
  
  Let us denote $\Omega_2 = \Omega_1 \cap \cbr{\omega \in \Omega: [X]_T\rbr{\omega} <  \infty}$. Naturally, $\P \rbr{\Omega_2} =1$. For $\omega \in \Omega_2$ we also have $\sum_{0 < s \le t} \rbr{\Delta X_s(\omega)}^2 <\infty$. This observation together with Remark~\ref{condition_jumps} yields that if $\omega \in \Omega_2$ and $x = X(\omega)$ then the sequence $\rbr{J_t (x^n, \cdot)}_n$ converges weakly in $L^1\rbr{\R}$ to $J_t (x, \cdot)$. 

  Thus we proved that for $\omega \in \Omega_2$ and $x = X\rbr{\omega}$ both required (weak) convergences hold, thus $x \in \mathbb{L}_1^S\rbr{\rbr{c_n}_n}$.
\end{proof}

We now begin the proof of Theorem~\ref{conv}. It is achieved via several lemmas. \emph{From now on, $X=(X_t)_{t\in [0,\infty)}$ will be a c{\`a}dl{\`a}g semimartingale, and $\mathcal{L}_t^u$, $t\geq 0$, $u\in \R$, its local times}. We will also need to consider, given $d>0$ and $z\in\R$, the semimartingale $X^{z,d}$, the processes $Y^{z,d},\tilde{X}^{z,d}$,  the functions $F^{z,d}, I^d_z\colon \R \to \R$ and the sequence of stopping times $(\tau_{n}^{z,d})_{n\in \N}$  defined as follows:
\begin{align*}
  &X^{z,d}_t:=F^{z,d}\rbr{X_t},   \quad \text{where}  \quad    F^{z,d}(x) := (z-d/2) \vee \left(x \wedge \left(z+d/2\right)\right) , \\
  &Y^{z,d}_t :=  X^{z,d}_t - X^{z,d}_{t-} - I^d_z\rbr{X_{t-}} \Delta X_t ,  \quad  \text{ where } \quad    I^d_z(x):= \1_{(z-d/2,z+d/2]}(x)  , \\
  &\tilde{X}^{z,d}_t:=\sum_{n=1}^{\infty}X_{\tau_{n-1}^{z,d}}^{z,d}{\bf 1}_{\left[\tau_{n-1}^{z,d},\tau_{n}^{z,d}\right)}\left(t\right),
\end{align*}
where the sequence of stopping times is defined by induction as follows: $\tau_{0}^{z,d}:=0$,
\begin{align*}  
  \tau_{1}^{z,d} & :=\begin{cases}
  \inf\left\{ s>0:X_{s}^{z,d}\in\left\{ z-d/2,z+d/2\right\} \right\}  & \text{ if }X_{0}^{z,d}\notin\left\{ z-d/2,z+d/2\right\} \\
  \inf\left\{ s>0:\left|X_{s}^{z,d}-X_{0}^{z,d}\right|=d\right\}  & \text{ if }X_{0}^{z,d}\in\left\{ z-d/2,z+d/2\right\} 
  \end{cases},
\end{align*}
and, for $n\geq 1$, 
\begin{equation*}
  \tau_{n+1}^{z,d}:=\begin{cases} 
  \inf\left\{ s>\tau_{n}^{z,d}:\left|X_{s}^{z,d}-X_{\tau_{n}^{z,d}}^{z,d}\right|=d\right\}  & \text{ if }\tau_{n}^{z,d}<\infty\\
  \infty & \text{ if }\tau_{n}^{z,d}=\infty
  \end{cases},
\end{equation*}
where we apply the usual conventions $\inf\emptyset:= \infty$ and $[\infty,\infty) := \emptyset$.
  
The first step in the proof of Theorem~\ref{conv} is to obtain a convenient formula for the quantity to be estimated, as we will now do.

\begin{lemma}\label{thm: long formula for difference}
  There exists a c{\`a}d{\`a}g adapted process $R^{z,d}$ with values in $(-2,0]$ such that
  \begin{align*}
     Q_t^{z,d} = R_{t}^{z,d} d + \frac{1}{d}\sum_{0<s\le t}\left(\Delta X_{s}^{z,d}\right)^{2} + \frac{2}{d} \int_{0}^{t}\left(X_{s-}^{z,d}-\tilde{X}_{s-}^{z,d}\right)\dd{X^{z,d}_s},  \quad t\in [0,\infty).
  \end{align*}
\end{lemma}

In the proof of Lemma~\ref{thm: long formula for difference}, and later on, we will make use of the following simple lemma.

\begin{lemma}\label{thm: Tanaka to Fzd(X)}
  $ X^{z,d}$ is a semimartingale and the following identity holds
  \begin{align*}
    X^{z,d}_t -X^{z,d}_0 = &  \frac{1}{2}\left( {\mathcal L}_t^{z-d/2} - {\mathcal L}_t^{z+d/2} \right) +  \int_0^t I^d_z\rbr{X_{s-}} \dd X_s  + \sum_{0<s\le t}  Y^{z,d}_s,
    \quad t\in [0,\infty). 
  \end{align*}
\end{lemma}

\begin{proof}
  Expressing $F^{z,d}$ as
  \begin{align}\label{eq: Fzd}
    F^{z,d}(x) = z + \frac{1}{2}\left|x-\rbr{z-d/2} \right| - \frac{1}{2}\left|x-\rbr{z+d/2} \right|
  \end{align}
  shows that it equals the difference of convex functions, and allows to quickly calculate its  derivatives as follows. Its first left-derivative  is $\rbr{F^{z,d}}'(x)=I^d_z(x)$,  and its  second (distributional) derivative is $\rbr{F^{z,d}}''=\delta_{z-d/2} -\delta_{z+d/2}$. The thesis  follows from the Tanaka--Meyer formula applied to $X$ and $F^{z,d}$ (see e.g. \cite[Chapter~IV, Theorem~70]{Protter2004}).
\end{proof} 

\begin{proof}[Proof of Lemma~\ref{thm: long formula for difference}]
  We have 
  \begin{align}\label{eq: crossing X same as Xzd}
    \cross{z,d}X{\left[0,t\right]}=\cross{z,d}{X^{z,d}}{\left[0,t\right]}.
  \end{align}
  We can now define 
  \begin{equation}\label{eq:jed}
    R_{t}^{z,d}:=  \cross{z,d}{X^{z,d}}{\left[0,t\right]}-\frac{1}{d^{2}}\sum_{n=1}^{\infty}\left(X_{\tau_{n}^{z,d}\wedge t}^{z,d}-X_{\tau_{n-1}^{z,d}\wedge t}^{z,d}\right)^{2} .
  \end{equation}
  Notice that $R_{t}^{z,d} \in\left(-2, 0\right]$ since only the first and last non-zero term in the above sum may differ from $d^2$, and they are then strictly smaller than $d^2$. Let us now work out an alternative expression for $ R_{t}^{z,d}$. Using integration by parts we get 
  \begin{equation}\label{eq: discrete [X] minus [X]}
    \sum_{n=1}^{\infty}\left(X_{\tau_{n}^{z,d}\wedge t}^{z,d}-X_{\tau_{n-1}^{z,d}\wedge t}^{z,d}\right)^{2}-\left[X^{z,d}\right]_{t}=2\int_{0}^{t}\left(X_{s-}^{z,d}-\tilde{X}_{s-}^{z,d}\right)\dd{X^{z,d}_s},
  \end{equation}
  where $\left[X^{z,d}\right]$ denotes the quadratic variation of $X^{z,d}$. Lemma~\ref{thm: Tanaka to Fzd(X)} shows that 
  $$   
    X^{z,d}_t - X^{z,d}_0 - \int_0^t I^d_z\rbr{X_{s-}} \dd X_s 
  $$
  is a process of finite variation and thus, denoting by $\left[Y\right]^{c}$ the continuous part of the quadratic variation of the semimartingale $Y$, we get that 
  \begin{align*}
    \left[X^{z,d}\right]_{t}^{c}
    =\sbr{\int_0^\cdot I^d_z(X_{s-}) \dd X_s }_t^{c}
    = \int_0^t I^d_z(X_{s-}) \dd [X]^{c}_s ,
  \end{align*} 
  and so by the occupation formula \eqref{eq:odf} we get 
  \begin{align}
    \left[X^{z,d}\right]_{t}  =\left[X^{z,d}\right]_{t}^{c}+\sum_{0<s\le t}\left(\Delta X_{s}^{z,d}\right)^{2} =\int_{z-d/2}^{z+d/2}{\mathcal L}_{t}^{u}\dd u+\sum_{0<s\le t}\left(\Delta X_{s}^{z,d}\right)^{2}.\label{eq:trz}
  \end{align}
  Combining \eqref{eq: crossing X same as Xzd}, \eqref{eq:jed}, \eqref{eq: discrete [X] minus [X]}, and \eqref{eq:trz} yields the thesis.
\end{proof}

To take advantage of the formula in Lemma~\ref{thm: long formula for difference}, we need a more convenient expression for the integral with respect to $X^{z,d}$; to obtain one, we again employ Lemma~\ref{thm: Tanaka to Fzd(X)}. This leads us to have to estimate the integral in $\dd z$ of three stochastic integrals (with respect to $\mathcal{L}^{z + d/2}-\mathcal{L}^{z - d/2}$, $X$ and $\sum_{0<s\le \cdot }  Y^{z,d}_s$, respectively); we will now do that, using a lemma for each integral.

\begin{lemma}\label{thm: local time no contribution}
  For $t\in [0,\infty)$, one has
  $$
    \int_0^t \left(X_{s-}^{z,d}-\tilde{X}_{s-}^{z,d}\right)\dd{\mathcal L}_s^{z\pm d/2} = 0.
  $$
\end{lemma}

\begin{proof}
  By \cite[Chapter~IV, Theorem~69]{Protter2004}), each of the atomless measures $\dd{\mathcal L}_s^{z\pm d/2}$ is carried by the corresponding set
  \begin{equation*}
    \cbr{s>0: X_s = X_{s-} = z\pm d/2} ,
  \end{equation*}
  and since the sets $$\cbr{s>0: X_s = X_{s-} = z\pm d/2 \neq \tilde{X}_{s-}^{z,d} }$$ are countable (because they are subsets of the jumps of the c{\`a}dl{\`a}g process $\tilde{X}^{z,d}$), we conclude that  $\dd{\mathcal L}_s^{z\pm d/2}$ is carried by the set
  $$
    \cbr{s>0: X_s = X_{s-} = z\pm d/2 = \tilde{X}_{s-}^{z,d} } \subseteq \cbr{s>0:  X^{z,d}_{s-} = \tilde{X}_{s-}^{z,d} } .
   $$ 
\end{proof} 

The stochastic integral with respect to $X$ will be estimated using the following lemma. 

\begin{lemma}\label{thm: bound for martingale integral}
  Let $(H^z)_{z\in \R}$ be a measurable\footnote{We mean that the function $(z,\omega,t)\mapsto H^z_t(\omega)$ is $\mathcal{B}(\R)\times \mathcal{P}$-measurable, where $\mathcal{B}(\R)$ are the Borel sets and $\mathcal{P}$ the predictable $\sigma$-algebra.} family of predictable process and assume that there exist constants $d,M\in (0,\infty)$   s.t., for all $s\geq 0$, $ |H^z_s|\leq d$ a.s. for all $z\in \R$, and  $H_s^z= 0$ a.s.  for all $|z| > M + d/2$. Given a semimartingale $S$, define
  $$
    W_t(S):=W_t^{z,d}(S):= \int_{\R}\left| \frac{2}{d} \int_{0}^{t}H^z_s I^d_z\rbr{X_{s-}} \dd S_s \right| \dd z.
  $$  
  If we assume that $S=V$ has finite variation $|V|_t<\infty$ a.s., then 
  \begin{align*}
    |W_t(V)| \leq 2 d \cdot |V|_t ,
  \end{align*}
  whereas if $S=N$ is a martingale s.t. $\E [N]^p_t<\infty$, then there exists a constant $C_p\in (0,\infty)$ (which depends only on $p$) s.t.
  \begin{align}\label{eq: bound for int mart}
    \E (W_t(N))^{2p} \leq \left(2M+d\right)^{2p-1}d C_p  \E (N)_t^p .  
  \end{align}
\end{lemma}

\begin{proof}
  Since $I^d_z\rbr{X_{s-}}=1$ if $X_{s-} - \frac{d}{2} \leq z <  X_{s-} + \frac{d}{2}$, and $I^d_z\rbr{X_{s-}}=0$ otherwise, we get 
  \begin{align}\label{eq: int Iz}
    \int_{\R} I^d_z\rbr{X_{s-}}  \dd z=d
  \end{align} 
  and so by Fubini's theorem
  \begin{align*}
    |W_t(V)| \leq 2 \int_{\R} \rbr{ \int_{0}^{t} I^d_z\rbr{X_{s-}} \dd \left| V \right|_s} \dd z  = 2d \left| V \right|_t .
  \end{align*}
  Since $(\int  |g| \dd \mu)^p \leq \int |g|^p \dd \mu $ holds for any probability measure $\mu$, we get that 
  \begin{align}\label{eq: p inside}
    \left(\int_\Omega  |g| \dd \mu\right)^p \leq \mu(\Omega)^{p-1}\int_\Omega |g|^p \dd \mu 
  \end{align} 
  for any positive finite measure $\mu$ on $\Omega$. Since
  $$  
    \textstyle G^z_t:= \frac{2}{d} \int_{0}^{t}H_s^{z} I^d_z \rbr{X_{s-}} \dd N_s
  $$
  equals $0$ when $|z| > M + d/2$, applying \eqref{eq: p inside} it follows that 
  \begin{align}\label{eq: 1st ineq mart bound}
    \E (W_t(N))^{2p} =
    \E \left(\int_{-M - d/2}^{M + d/2}  |G_t^z| \dd z \right)^{2p} 
    \leq (2M + d)^{2p-1} \E\left ( \int_{\R}  \left( G_t^z  \right)^{2p}\dd z\right) .
  \end{align} 
  Burkholder-David-Gundy inequality applied to $ G^z$ gives that 
  \begin{align}\label{krokkk}
    \E \left( G_t^z  \right)^{2p}
    \leq c_p \E \left( \sqrt{[G^z]_t}  \right)^{2p}
    = c_p\E \left ( \int_{0}^{t} \Big(\frac{2}{d} H_s^{z} I^d_z \rbr{X_{s-}} \Big)^2  \dd [N]_s  \right)^{p} 
    =:A_t^z.
  \end{align}
  Using first \eqref{eq: p inside}, and then $|H|\leq d$ and $(I^d_z)^{2p}= I^d_z$, we get the two inequalities
  \begin{align}\label{eq: A}
    A_t^z \leq c_p \E   \left ( [N]_t^{p-1} \int_{0}^{t} \Big(\frac{2}{d} H_s^{z} I^d_z \rbr{X_{s-}} \Big)^{2p}  \dd [N]_s  \right)  \leq 4^p c_p \E  \Big( [N]_t^{p-1}  \int_{0}^{t} I^d_z \rbr{X_{s-}}   \dd [N]_s \Big) .
  \end{align} 
  Applying Fubini's theorem and combining \eqref{eq: 1st ineq mart bound}, \eqref{krokkk} and \eqref{eq: A} we get
  \begin{align*}
    \E (W_t(N))^{2p}  \leq 4^p c_p (2M + d)^{2p-1} \E  \Big( [N]_t^{p-1}  \int_{0}^{t}\big ( \int_{\R}  I^d_z \rbr{X_{s-}}  \dd z \big) \dd [N]_s \Big) ,
  \end{align*} 
  and now \eqref{eq: int Iz} yields the thesis with $C_p:= 4^p c_p$.
\end{proof} 

To deal with the stochastic integral with respect to $\sum_{0<s\le \cdot } Y^{z,d}_s $ we will use the following lemma.

\begin{lemma}\label{le: int wrt Y}
  $$ 
    \int_\R \dd z \sum_{0<s\leq t}  |Y^{z,d}_s|
    \leq \sum_{0<s\leq t} (\Delta X_s)^2 g_d(\Delta X_s), 
  $$
  where 
  \begin{align*}
    g_d(x):= \begin{cases}
    1 & \text{ if }   |x|\le d  \\
    \frac{5d}{|x|} & \text{ if }  |x| > d 
    \end{cases},
    \quad x\in \R , d>0 .
  \end{align*}
\end{lemma} 

\begin{proof}
  Since the expression $\Delta f(X_{s})-f^\prime (X_{s-} )\Delta X_s$ is linear in $f$ and equals 0 when $f$ is a constant, using \eqref{eq: Fzd} and \eqref{eq:Two expressions for J} shows that
  \begin{align*}
    Y^{z,d}_s = Z_s^{z, -d} - Z_s^{z, d} , \text{ where }  Z_s^{z,u}:= \left|X_s - \left(z + u/2\right)\right| \1_{  \llbracket X_{s-}, X_{s} \rrparenthesis  } \left(z +u/2 \right), \quad z,u\in \R, 
  \end{align*}
  and to conclude we only need to compute the $L^1$-norm 
  \begin{align*}
    \sum_{0<s\leq t}  \int_\R   \dd z  \,| Z_s^{z, - d} - Z_s^{z, d}|
  \end{align*}
  of $Z_s^{z, - d} - Z_s^{z, d}$. 
  
  If $d \ge |\Delta X_s|$, then \eqref{eq:int |b-u|^p} with $p=1$ gives that
  \begin{align}\label{eq: L1 norm Z-Z small jump}
    \int_\R   \dd z  |Z_s^{z, - d} - Z_s^{z, d}| \le \int_\R   \dd z  \rbr{ |Z_s^{z,-  d}| + | Z_s^{z, d}|} = (\Delta X_s)^2 .
  \end{align} 
  
  To deal with the case $d < |\Delta X_s|$, let us notice that 
  \begin{align*} 
    \text{ if } z +d/2 < X_{s-}\wedge X_s  \text{ or } z -d/2 > X_{s-}\vee X_s \text{ then } Z_s^{z, - d} = Z_s^{z, d} = 0  
  \end{align*} 
  and 
  \begin{align*} 
    \text{ if }  X_{s-}\wedge X_s < z -d/2 < z +d/2 < X_{s-}\vee X_s \text{ then } |Z_s^{z, - d} - Z_s^{z, d}| \le d. 
  \end{align*} 
  The last estimate follows from the fact that if 
  $$
    X_{s-}\wedge X_s < z -d/2 < z +d/2 < X_{s-}\vee X_s,
  $$
  then $Z_s^{z,-d}= \left|X_s - \left(z - d/2\right)\right|$,   $Z_s^{z,d}= \left|X_s - \left(z + d/2\right)\right|$, and the inequality $| |a| - |b| | \le |a-b|$. 

  Finally, in the case 
  $$
    z \in \sbr{X_{s-}\wedge X_s -d/2 , X_{s-}\wedge X_s + d/2 } \cup \sbr{X_{s-}\vee X_s -d/2 , X_{s-}\vee X_s + d/2 }  
  $$
  we apply the estimate $|Z_s^{z,u}| \le |\Delta X_s| $, valid for any $z,u \in \R$.
 
  Putting together three considered cases we have the estimate 
  \begin{align}\label{eq: L1 norm Z-Z big jump}
    \begin{split}
      \int_\R \dd z  |Z_s^{z, - d} - Z_s^{z, d}|  \leq  & \int_{X_{s-}\wedge X_s - d/2 }^{X_{s-}\wedge X_s + d/2}\dd z \rbr{ |Z_s^{z,-  d}| + | Z_s^{z, d}|} + \int_{X_{s-}\wedge X_s + d/2 }^{X_{s-}\vee X_s - d/2} d \dd z  \\
      &+ \int_{X_{s-}\vee X_s - d/2 }^{X_{s-}\vee X_s + d/2} \dd z \rbr{ |Z_s^{z,-  d}| + | Z_s^{z, d}|} \\
     \le  & 2 d  |\Delta X_s| + \rbr{ |\Delta X_s| - d}d+ 2 d  |\Delta X_s| \le 5 d  |\Delta X_s|.
    \end{split}
  \end{align} 
  From \eqref{eq: L1 norm Z-Z small jump}, \eqref{eq: L1 norm Z-Z big jump} it follows that the $L^1$-norm of $Z_s^{z, - d} - Z_s^{z,d}$ is bounded by 
  \begin{align*}
    \sum_{0<s\leq t: d \ge  |\Delta X_s|} (\Delta X_s)^2  + 
    \sum_{0<s\leq t }  \1_{(0,|\Delta X_s |)}\left(d\right) 5  d |\Delta X_s|  = \sum_{0<s\leq t} (\Delta X_s)^2 g_d(\Delta X_s) , 
  \end{align*} 
  which concludes the proof.
\end{proof} 

\begin{proof}[Proof of Theorem~\ref{conv}]
  For now assume that $X$ is in $\cS^{2p}$ for some $p\in [1,\infty)$, and $|X|$ is bounded by a constant $M$. Since $\cross{z,d}X{\left[0,t\right]}$ and ${\mathcal L}_{t}^{u}$ are equal $0$ for $\left|z\right|,\left|u\right|>M_t:=\sup_{0\le s\le t}\left|X_{s}\right|<\infty$, $Q_t^{z,d}=0$ for any $z \notin [-M_t -d/2, M_t +d/2]$, and thus 
  $$
    \int_\R \left| Q^{z,d}_t \right| \dd z
    = \int_{-M_t -d/2}^{M_t +d/2} \left| Q^{z,d}_t \right|\dd z .
  $$
  We can now apply Lemma~\ref{thm: long formula for difference} to estimate the latter as a sum of three terms. The first is
  \begin{equation}\label{eq:estim0}
   \left| \int_{-M_t -d/2}^{M_t +d/2} d \cdot  R_{t}^{z,d} \dd z  \right|
    \le 2d\cdot 2 \cdot \left(M_t +d/2\right) 
    = 4M_td+2 d^2 .
  \end{equation}
  The second term is  
  \begin{align}\label{eq: dx sum [X]}
    \int_{-M_t -d/2}^{M_t +d/2}  \frac{1}{d} \sum_{0<s\le t}\left(\Delta X_{s}^{z,d}\right)^{2} \dd z
    = \sum_{0<s\le t} \int_{-M_t -d/2}^{M_t +d/2}  \frac{1}{d} \left(\Delta X_{s}^{z,d}\right)^{2} \dd z .
  \end{align} 
  By definition of $X_{s}^{z,d}$ we have that $|\Delta X_{s}^{z,d}|\leq |\Delta X_s|\wedge d$, and if $X_{s-}<X_s$ then $\Delta X_{s}^{z,d}=0$ whenever $z \notin [X_{s-} - \frac{d}{2}, X_s+  \frac{d}{2}]$, and analogously if  $X_{s}\leq X_{s-}$ then $\Delta X_{s}^{z,d}=0$ whenever $z \notin [X_{s} - \frac{d}{2}, X_{s-} +  \frac{d}{2}]$. This gives the first of the following inequalities  
  \begin{align}\label{eq: bound dx sum [X]}
    \int_{\R}   \frac{1}{d} \left(\Delta X_{s}^{z,d}\right)^{2}  \dd z
    \leq \frac{1}{d}  \Big( |\Delta X_{s}| \wedge d \Big)^2 \Big( |\Delta X_{s}| +  d  \Big)
    \leq \frac{1}{d}  \Big( |\Delta X_{s}| \wedge d \Big)^2 2( |\Delta X_{s}| \vee  d) .
  \end{align} 
  Thus, using the identity
  \begin{align*}
    \frac{1}{d} \Big( |\Delta X_{s}| \wedge d \Big)^2 \Big(|\Delta X_{s}| \vee  d \Big) = (\Delta X_{s})^2 \wedge  (d |\Delta X_{s}|) ,
  \end{align*} 
  which can easily verified separately for the cases $d< |\Delta X_{s}|$ and $d\geq |\Delta X_{s}|$, combined with \eqref{eq: bound dx sum [X]} and \eqref{eq: dx sum [X]}, gives that
  \begin{align}\label{eq:estim1}
    \int_{-M_t -d/2}^{M_t +d/2} \frac{1}{d} \sum_{0<s\le t} \left(\Delta X_{s}^{z,d}\right)^{2} \dd z
    \leq 2 \sum_{0<s\le t} (\Delta X_{s})^2 \wedge  (d |\Delta X_{s}|)  =: D^d_t.
  \end{align} 

  The third and last term which we need to estimate is
  \begin{align}\label{eq: stoch integral to estimate hard}
    \int_\R\left| \frac{2}{d} \int_{0}^{t}\left(X_{s-}^{z,d}-\tilde{X}_{s-}^{z,d}\right)\dd{X^{z,d}_s} \right|  \dd z  .
  \end{align} 
  To do so, we use Lemma~\ref{thm: Tanaka to Fzd(X)}  to write this as the sum of the  integrals with respect to $\mathcal{L}^{z + d/2}-\mathcal{L}^{z - d/2}$, $X$, and $\sum_{0<s\le \cdot } Y^{z,d}_s$. The first  integral is zero, thanks to Lemma~\ref{thm: local time no contribution}. To estimate the second integral (in $\dd X$), we write  the canonical semimartingale decomposition $X=N+V$ of $X$ as a local martingale $N$ and a predictable process of finite variation $V$, and apply Lemma~\ref{thm: bound for martingale integral} with 
  $$
    H^z_s:=X_{s-}^{z,d}-\tilde{X}_{s-}^{z,d}.
  $$ 
  Assumptions of Lemma~\ref{thm: bound for martingale integral} are satisfied because
  \begin{align}\label{eq: diff X X tilde}
    \left|X_{s-}^{z,d}-\tilde{X}_{s-}^{z,d} \right| \leq d  \quad \text{for all} d,s>0,z\in \R,
  \end{align}
  $X\in \cS^{2p}$ implies $[N]_\infty \in L^p\rbr{\R}, |V|_\infty \in L^{2p}\rbr{\R}$, and from the implication
  \begin{align*}
    \text{if}\quad X_{s-} \notin \sbr{z - d/2, z +d/2} & \quad \text{then}  \quad X_{s-}^{z,d}-\tilde{X}_{s-}^{z,d} = 0
  \end{align*}
  it follows that $X_{s-}^{z,d}-\tilde{X}_{s-}^{z,d} = 0$ unless $|z -X_{s-}|\leq d/2$, and since the constant $M$ satisfies  $ M \ge\sup_{0\le s}\left|X_{s}\right| $, this implies that   $X_{s-}^{z,d}-\tilde{X}_{s-}^{z,d} = 0$  for $|z| > M + d/2$. To estimate the third integral, we apply Lemma~\ref{le: int wrt Y} and \eqref{eq: diff X X tilde}. Combining these three estimates we can bound the third term (the one in \eqref{eq: stoch integral to estimate hard}), and this bound, combined with those obtained in \eqref{eq:estim0} and \eqref{eq:estim1} for the first and second terms gives that
  \begin{align}\label{eq: combo of all estimates}
    \int_\R |Q_t^{z,d}| \dd z 
    \leq 4M_t d+2d^2+ D^d_t + 2 d |V|_t + W_t^{z,d}(N) + 2\sum_{0<s\leq t} (\Delta X_s)^2 g_d(\Delta X_s)  .
  \end{align} 
  Clearly \eqref{eq: conv S2} follows from \eqref{eq: combo of all estimates} and the dominated convergence theorem, since any $X$  in $\cS^{2p}$ satisfies $\E([X]_t^p+ \sup_{s\leq t} |X_s|^{2p}  )<\infty$ (see e.g. \cite[Chapter~7, Section~3, Number~98, Page~295, Equations~98.5 and~98.7]{Dellacherie1982}), $g_d$ satisfies $0\leq g_d(x)\leq 5$ and $g_d(x)\to 0$ at all $x\neq 0$ as $d \downarrow 0$, and we use the fact that $\E  (W_t^{z,d}(N))^{p}\leq \sqrt{ \E  (W_t^{z,d}(N))^{2p}}$, which goes to zero thanks to Lemma~\ref{thm: bound for martingale integral}.
 
  If $X$ is an arbitrary semimartingale, to prove \eqref{eq:conv} we can assume w.l.o.g. that   $X$ is in $\cS^2$ by pre-localisation, see e.g. \cite[Chapter~IV, Theorem~13]{Protter2004}). Now let $d_k\geq 0$ be such that $\sum_{k=1}^{\infty} d_k < \infty$; the term $ W_t^{z,d}(N)$ is controlled by the estimate~\eqref{eq: bound for int mart}, which gives that 
  $$
    \E\left ( \sum_{k=1}^\infty \rbr{ W_t^{z,d_k}(N)}^{2}\right)<\infty,
  $$
  from which we conclude that $W_t^{z,d_k}(N)\to 0$ a.s. as $k \to \infty$. As the remaining terms in \eqref{eq: combo of all estimates} converge to $0$ a.s. as $k\to \infty$ (the term $D^{d_k}_t$ and the last term by dominated convergence, the others trivially), \eqref{eq:conv} follows.  
\end{proof}

\bibliography{references}{}
\bibliographystyle{amsalpha}

\end{document}